\documentclass[12pt]{amsart}
\usepackage{amscd, amssymb,latexsym,amsmath, amscd, amsmath}
\usepackage[all]{xy}
\usepackage{anysize}
\usepackage[sc]{mathpazo} 
\usepackage{color}
\marginsize{2.5cm}{2.5cm}{2.5cm}{2.5cm}
\newtheorem{theorem}{Theorem}[section] 
\newtheorem{lemma}[theorem]{Lemma}
\newtheorem{corollary}[theorem]{Corollary}
\newtheorem{proposition}[theorem]{Proposition}
\theoremstyle{definition}
\newtheorem{conjecture}[theorem]{Conjecture}

\theoremstyle{remark}
\newtheorem{remark}{Remark}

\begin{document}

\title[Newforms as test vectors for representations of Saito-Kurokawa type]{Newforms as test vectors for representations of Saito-Kurokawa type}

\author{U. K. Anandavardhanan}

\address{Department of Mathematics, Indian Institute of Technology Bombay, Mumbai - 400076, India.} \email{anand@math.iitb.ac.in}

\subjclass{Primary 20C15; Secondary 20C20, 20C33}

\date{}

\begin{abstract}
We study irreducible (admissible) representations $\pi$ of the split group $G_n=SO(2n+1)$, both over finite and $p$-adic fields, which are distinguished for the (unramified) nonsplit form of $K_n=SO(2n)$ and for one more subgroup $H_n$. Over finite fields, $H_n$ is the split $SO(2n)$ while over $p$-adic fields it is the Gross-Tsai paramodular subgroup (of minimal level for $\pi$). These are all multiplicity one situations. With $0 \neq \ell_{n,K} \in {\rm Hom}_{K_n}(\pi_n,1)$ and $0 \neq v_{n,H} \in \pi_n^{H_n}$, both being one dimensional spaces, we are interested in the vanishing/nonvanishing of $|\ell_{n,K}(v_{n,H})|^2$. Appropriately normalized, this defines the notion of the correlation coefficient $c(\pi;H_n,K_n)$ of $H_n$ and $K_n$ in $G_n$ with respect to $\pi$. Via the theta correspondence between ${\rm SL}(2)$ (or ${\rm Mp}(2)$), we establish a bijection $\pi_n \leftrightarrow \pi_1$ between representations of $G_n$ and $G_1$ under consideration. Our main results then prove that $d_n c(\pi_n;H_n,K_n)$ is independent of the rank $n$, where $d_n =\dim \pi_n$ for finite fields and is defined in terms of the index of $H_n$ in the standard maximal compact subgroup of $G_n$ in the $p$-adic case. We also highlight a conjecture from \cite{ana25} on Legendre polynomials modulo $p$. Its truth will have the arithmetic consequence that over $\mathbb F_p$, given $\pi_1$, we can choose $u \in \mathbb F_p$ such that, in every rank $n$, all the correlation coefficients $c(\pi_n;H_n,u_nK_n u_n^{-1})$ are simultaneously nonzero and moreover they belong to $\mathbb Z_p^\times$ under any identification $\mathbb C \simeq \overline{\mathbb Q}_p$, where $u_n$ is a specific unipotent element in $G_n$ corresponding to $u$.
\end{abstract}

\maketitle

\section{Introduction}\label{intro}

For a group $G$ and two subgroups $H$ and $K$, roughly speaking, we are interested in the question of correlating the phenomenon of distinction, in the relative Langlands programme, for an irreducible representation of $G$ with respect to the subgroups $H$ and $K$. This is especially interesting when the representation admits a one dimensional space of invariant linear forms for each of the two subgroups. In this paper, $G$ is $G_n = SO(n+1,n)(F)$, the split odd special orthogonal group, where $F$ is either a finite field or a $p$-adic field of odd (residue) characteristic. By a representation of Saito–Kurokawa type we mean a representation $\pi$ of $G_n$ which is a nonzero small theta lift of an irreducible (resp. irreducible tempered genuine) representation of ${\rm SL}_2(F)$ (resp. ${\rm Mp}_2(F)$). For $n=2$, these are the so called Saito–Kurokawa representations. For $K$, we take $K_n=SO(n+1,n-1)(F)$, where $SO(n+1,n-1)$ denotes the quasi-split, non-split even special orthogonal group of absolute rank $n$ whose splitting field is quadratic over $F$ and is also unramified when $F$ is $p$-adic. Writing $H_n=SO(n,n)(F)$ in both the cases, when $F$ is finite, we take $H=H_n$, and when $F$ is $p$-adic, we take $H$ to be a specific conjugate of the Gross–Tsai paramodular subgroup $P_{n,m}$ at a level $m$ determined by $\pi$. Thus, we are investigating correlation in a GGP context. The results of this paper can be summarised as  
{\it correlation at a higher rank is completely captured at rank one}
in the cases that we investigate.

The present work is motivated by two results which we first recall. The first one is stated in \cite[\S 7.1]{bhko24} where it is shown that there is an embedding of the unramified torus in ${\rm PGL}_2(\mathbb Q_p)$, depending only on the conductor of the representation, such that the local Atkin-Lehner newform is a test vector for the corresponding local unramified toric period, for all irreducible supercuspidal representations of ${\rm PGL}_2(\mathbb Q_p)$ of a given conductor. This is done for conductor (which is necessarily even for the representation to have a nonvanishing local unramified toric period) $\geq 4$. The missing case of conductor two can be approached, as explained in \cite[pp. 85-86]{vat23}, by reformulating the problem as one of simultaneous nonvanishing of the correlation coefficient \cite[\S 8, \S9]{gro91}
\[c(\pi;H,K) = \frac{|(v_H,v_K)|^2}{(v_H,v_H)(v_K,v_K)}\]
for all irreducible cuspidal representations $\pi$ of trivial central character of ${\rm GL}_2(\mathbb F_p)$ for a specific embedding $K$ of $\mathbb F_{p^2}^\times$ in ${\rm GL}_2(\mathbb F_p)$ and $H$ being the standard diagonal embedding of the split torus; here $v_K$ and $v_H$ are the respective fixed vectors and $(\cdot,\cdot)$ is an invariant Hermitian inner product on $\pi$. Following a strategy initially implemented in \cite{aj22}, this is carried out in \cite{ana25} where it is shown that the answer is in terms of mod $p$ reduction of Legendre polynomials. In all these works $p$ is assumed to be odd.

The higher depth case treated in \cite[\S 7]{bhko24} requires an intricate analysis on the Kirillov model and this is done for $Q_p$. In \cite{ana25}, the depth zero case is taken up for all $p$-adic fields $F$, where the crucial input is the mod $p$ reduction of irreducible complex representations of ${\rm GL}_2(\mathbb F_q)$, where $q$ is the residue cardinality of $F$, and this is due to Diamond \cite{dia07}. 

The key observation leading to the present work is that the main result in \cite{ana25} does not need the representation of ${\rm PGL}_2(\mathbb F_q)$ to be cuspidal. The result there is uniformly stated 
for all irreducible representations of ${\rm PGL}_2(\mathbb F_q)$ which have nontrivial fixed vectors for both the split and nonsplit tori. Their number is $q$, the cardinality of the corresponding double coset space. Since $|SO^+(2n,\mathbb F_q) \backslash SO(2n+1,\mathbb F_q)/SO^-(2n,\mathbb F_q)| = q$, it is natural to try to compute the correlation coefficient for this triple. It is not difficult to identify the $q$ many irreducible representations of $G_n = SO(2n+1,\mathbb F_q)$ which are distinguished for both $H_n = SO^+(2n,\mathbb F_q)$ and $K_n = SO^-(2n,\mathbb F_q)$; these turn out to be the big theta lifts from ${\rm SL}_2(\mathbb F_q)$ of precisely the $(q-1)$ many nontrivial irreducible representations with singleton $L$-packets (these big theta lifts are known to be irreducible, cf. \cite{kri24,kri25}), and of course the trivial representation of $G_n$. 

We initially thought perhaps the mod $p$ reduction strategy of \cite{ana25} can be attempted since spherical harmonics in terms of Gegenbauer polynomials as in \cite{bry92} provide a setup analogous to that of \cite{ana25}. Of course, unlike the case of $n=1$, the subgroup sizes are not coprime to $p$ when $n \geq 2$, and this creates a hurdle to mimic \cite{ana25}, where at the very start there is an averaging technique. In any case, experimenting with low rank character tables, it turned out that this strategy cannot be salvaged in the higher rank case. So we were forced to look at the question purely on the complex side without reducing mod $p$, employing the knowledge that we are looking only at theta lifts, and with a simple enough computation on the Schr\"odinger model, we get the answer; the nontrivial irreducible representations of $G_n$ and $G_1$ under consideration are in natural bijection $\pi_n\leftrightarrow\pi_1$ via the intermediary of ${\rm SL}_2(\mathbb F_q)$ by theta correspondence and then the correlation coefficient for the triple $(G_n,H_n,K_n)$ for an irreducible representation $\pi$ of $G_n$ is essentially the correlation coefficient for the triple $(G_1,H_1,K_1)$ for its image under this natural bijection. The case $n=1$ is answered in \cite{ana25}. In order to state it more precisely, as in \cite{ana25}, we need to consider a family of non-split groups $K_{n,\alpha,u}$ indexed by $u \in \mathbb F_q$ and $\alpha \in \mathbb F_q^\times \setminus \mathbb F_q^{\times 2}$ (see equation (\ref{k-def}) in \S \ref{ortho-para} for the precise definition). We also extend the above bijection $\pi_n \leftrightarrow \pi_1$, defined for nontrivial irreducible representations, by $1_{G_n}\leftrightarrow1_{G_1}$.

\begin{theorem}\label{thm1} 
Let $\pi_n$ and $\pi_1$ be irreducible representations of $G_n$ and $G_1$ in bijection as above.
Then, for every $u\in\mathbb F_q$, we have
\[
(\dim\pi_n) c(\pi_n;H_n,K_{n,\alpha,u})
=
(\dim\pi_1) c(\pi_1;H_1,K_{1,\alpha,u}).
\]
\end{theorem}

\begin{remark}\label{rmk-thm1}
For $\pi_1 \neq 1$, $\dim \pi_n = \frac{q^{2n}-1}{q^2-1}\dim \pi_1$ (see Proposition \ref{kriz} (1)).
\end{remark}

Having gotten such a clean looking identity over finite fields, we found it tempting to ask whether such an $n$ to $1$ reduction is possible also over $p$-adic fields, where the $n=1$ case is the work in \cite[\S 7.1]{bhko24}, and this is indeed the case. Thus, we have $G_n$, $H_n$, and $K_n$ as before but over a $p$-adic field $F$ of odd residue characteristic with ring of integers $\mathfrak o_F$, uniformizer $\varpi$, and residue cardinality $q$. Fix a nontrivial additive character $\psi$ of $F$ and $\alpha \in \mathfrak o_F^\times \setminus \mathfrak o_F^{\times 2}$, and let $E=F(\sqrt{\alpha})$. We write $K_{n,\alpha}$ for $K_n$; this is of absolute rank $n$ and Witt index $(n-1)$ and is split over the unramified quadratic extension $E$. Note that $K_{1,\alpha}\simeq E^\times/F^\times$.

Considering that every infinite dimensional irreducible admissible representation of ${\rm PGL}_2(F)$ has an Atkin-Lehner local newform, we ask: what are the irreducible admissible representations of $G_n$ which are distinguished with respect to $K_{n,\alpha}$ and with a fixed vector for the paramodular subgroup $P_{n,m}$ of $G_n$ of level $m$ \cite{tsa13,gro16}? 

In order to address this question we make use of theta correspondence. Write $\theta_{j,\psi}$ and $\theta_{j,\psi}^{\leftarrow}$ for the forward and backward small theta lifts between ${\rm Mp}_2(F)$ and $G_j=SO(j+1,j)(F)$. 

Using the results on theta correspondence \cite{gt16,ggp20} and on paramodular newforms and the structure of parabolically induced representations  \cite{jan96,rs07,che26}, together with a little bit of Mackey theory \cite[\S 15]{ggp12}, we show that every nontrivial representation $\pi$ satisfying the two distinction conditions above is of the form
$\pi \simeq \theta_{n,\psi}\bigl(\theta^{\leftarrow}_{1,\psi}(\tau)\bigr)$ 
for an irreducible admissible $K_{1,\alpha}$-distinguished representation $\tau$ of $G_1\simeq {\rm PGL}_2(F)$. Here $\theta^{\leftarrow}_{1,\psi}(\tau)$ is the Waldspurger correspondent of $\tau$ on ${\rm Mp}_2(F)$. When $n\ge 2$ and $\tau$ is tempered of positive conductor, $\pi$ can be realized as the Langlands quotient of a parabolically induced representation. This realization allows us to construct its paramodular newform recursively, starting with the Atkin–Lehner newform in $\tau$.

At this point, we state a multiplicity one result for these paramodular newforms. It is mostly essentially due to \cite{che26}; see the introductory words in \S \ref{section-cheng} before Proposition \ref{prop-cheng} for more details about how this result is arrived at. Its proof is in \S \ref{section-crsj}.

\begin{theorem}\label{thm-cheng}
Let $\sigma$ be an irreducible admissible genuine representation of ${\rm Mp}_2(F)$ such that
$\tau=\theta_{1,\psi}(\sigma)$
is nonzero and tempered.  Assume that its conductor exponent $c=a(\tau)$ is even. Then, $\pi_n=\theta_{n,\psi}(\sigma)\neq0$, $\pi_n^{P_{n,m}}=0$ for $0\le m<c$, and $\dim\pi_n^{P_{n,c}} = 1$. 
\end{theorem}

For $r\geq 0$, let $\mathcal P_{j,2r}=b_{j,r}P_{j,2r}b_{j,r}^{-1}$, where $b_{j,r}={\rm diag}(\varpi^r I_j,1,\varpi^{-r}I_j)$, and put $D_j(r)=q^{-jr}[G_j(\mathfrak o_F):\mathcal P_{j,2r}]$. Thus $D_j(0)=1$ and $D_j(r)=q^{j(r-1)}(q^j+1)$ for $r\geq 1$. At positive conductor $2r$, take $v_{1,2r}$ to be the translate by $b_{1,r}$ of an Atkin–Lehner newform. 
For $2\le j\le n$, normalize the newvector $v_{j,2r}$ fixed by the conjugated subgroup $\mathcal P_{j,2r}$ recursively by requiring the inducing section $F_j$ to satisfy $F_j(1)=v_{j-1,2r}$.
For $u \in \mathfrak o_F$, we consider a specific $K_{n,\alpha,u}$-invariant linear form on $\pi_n$ which we denote by $\Lambda_{n,\alpha,u}$ (see Remark \ref{rmk-k-def} for its definition). 
At rank one, $\Lambda_{1,\alpha,u}$ is the linear toric functional for the embedding considered in
\cite[\S7]{bhko24}. We normalize $\Lambda_{j,\alpha,u}$ recursively from the rank one toric period $\Lambda_{1,\alpha,u}$, using volume one Haar measures on the compact groups used in these constructions. 

Now we are in a position to state the $p$-adic analogue of Theorem \ref{thm1}. 

\begin{theorem}\label{thm2} 
Let $\pi_n$ be a nontrivial irreducible admissible representation of $G_n$. Suppose that $\pi_n$ is $K_{n,\alpha}$-distinguished and that $\pi_n^{P_{n,m}}\ne0$ for some $m\ge0$.
\begin{itemize}
\item[1.] There is a unique irreducible admissible representation $\tau$ of $G_1\simeq {\rm PGL}_2(F)$ such that
 $\pi_n\simeq\theta_{n,\psi}\bigl(\theta_{1,\psi}^{\leftarrow}(\tau)\bigr)$. 
Moreover, $\tau$ is $K_{1,\alpha}$-distinguished. Conversely, every $K_{1,\alpha}$-distinguished tempered representation $\tau$ of positive even conductor arises in this way from a unique $\pi_n$ satisfying the hypotheses above.
\item[2.] Suppose that $\tau$ is tempered and $a(\tau)=2r>0$. Then the minimal paramodular level of $\pi_n$ is $2r$ and its space of paramodular newvectors is one dimensional. Also, for every $u\in\mathfrak o_F$, we have
   \[
   D_n(r)\,\Lambda_{n,\alpha,u}(v_{n,2r})
   =
   D_1(r)\,\Lambda_{1,\alpha,u}(v_{1,2r}).
   \]
\end{itemize}   
\end{theorem}

Combining Theorem \ref{thm2} with the argument in the proof of \cite[Corollary 7.2(i)]{bhko24}, we obtain the following simultaneous nonvanishing result.

\begin{corollary}
Let $F=\mathbb Q_p$, with $p$ odd, and fix $\alpha\in\mathbb Z_p^\times \setminus \mathbb Z_p^{\times 2}$. Let $E=F(\sqrt\alpha)$. There exists $u \in\mathbb Z_p^\times$ depending only on $p$ and $\alpha$, with the following property. Let $\tau=\tau_\chi$ be a dihedral supercuspidal representation associated to a character $\chi$ of $E^\times$ satisfying
$\chi|_{F^\times}=\eta_{E/F}$, $a(\chi)=r\ge2$,
where $\eta_{E/F}$ is the quadratic character associated with $E/F$.
Then
\[
\Lambda_{n,\alpha,u}(v_{n,2r})\ne0.
\]
Thus the same parameter $u$ makes the paramodular newvector a test vector for every such $\chi$, every $r\ge2$, and every rank $n \geq 1$.
\end{corollary}

\begin{remark}
It is perhaps possible to have an analogue of Theorem \ref{thm2} when $\tau$ is an irreducible representation of $PD^\times$, where $D$ is the quaternion division algebra over $F$. If $\tau$ admits a nonzero unramified toric functional, say $\Lambda^\prime_{1,\alpha,0}$, its Jacquet–Langlands correspondent on ${\rm PGL}_2(F)$  must have odd conductor exponent. It is then known that there is a unique ramified torus $T$, up to conjugacy, for which $\tau$ has a nonzero fixed vector, say $v_{1,2r+1}$. The value of $\Lambda^\prime_{1,\alpha,0}(v_{1,2r+1})$ is computed in \cite{ap24}. The analogue then would be to recursively construct the linear forms $\Lambda^\prime_{n,\alpha,0}$  and the Bessel newforms $v_{n,2r+1}$ for $G_n$, in the sense of  \cite{oka21}, and then to compute $\Lambda^\prime_{n,\alpha,0}(v_{n,2r+1})$ in terms of $\Lambda^\prime_{1,\alpha,0}(v_{1,2r+1})$.
\end{remark}

Now we state the higher rank analogue of the argument in \cite[pp. 85-86]{vat23} which reduces the period value at the paramodular newform in the depth zero case to the corresponding correlation coefficient over finite fields. In the result below, by $u$ (resp. $\alpha$) we denote an element of $\mathfrak o_F$ (resp. $\mathfrak o_F^\times \setminus \mathfrak o_F^{\times 2}$) as well as its reduction to $\mathbb F_q$. Keep the hypotheses as in Theorem \ref{thm2} (2) and moreover assume $a(\tau)=2$ and thus $r=1$. Let $C_1(1)$ be the kernel of reduction $G_1(\mathfrak o_F) \to G_1(\mathbb F_q)$. Consider $\tau^{C_1(1)}$ as a representation of $G_1(\mathbb F_q)$, and let $\pi_1$ be its unique irreducible constituent with a nonzero $H_1(\mathbb F_q)$-fixed vector.
Let $\pi_n$ be the representation of $G_n(\mathbb F_q)$ corresponding to $\pi_1$ as in Theorem \ref{thm1}. 

\begin{theorem}\label{thm3}
Suppose $|\Lambda_{1,\alpha,u}(v_{1,2})|^2 = c(\pi_1;H_1,K_{1,\alpha,u})$ for every $u \in \mathfrak o_F$.
Then,
\[
D_n(1)^2\,\bigl|\Lambda_{n,\alpha,u}(v_{n,2})\bigr|^2
=
D_1(1)^2\,
\frac{\dim\pi_n}{\dim\pi_1}\,
c(\pi_n;H_n,K_{n,\alpha,u}).
\]
\end{theorem}

\begin{remark}
Note that the identity in Theorem \ref{thm3} is compatible with the identities in Theorems \ref{thm1} and \ref{thm2}.
\end{remark}

We now return to the finite field setting and combine Theorem \ref{thm1} with the rank one result on the mod $p$ reduction of the correlation coefficient in \cite{ana25}. For this, fix an identification of $\mathbb C$ with $\overline{\mathbb Q}_p$. Let $\chi$ be the character of $\mathbb F_q^\times$ whose reduction modulo $p$ is $x\mapsto x$, and let $\eta=\chi^{(q-1)/2}\circ\det$ be the unique quadratic character of ${\rm PGL}_2(\mathbb F_q)$. Choose a character $\nu$ of $\mathbb F_{q^2}^\times$ such that either $\nu$ or $\nu^q$ reduces modulo $p$ to $x\mapsto x$.
To each of the $q$ irreducible representations $\pi_1$ of $G_1\simeq {\rm PGL}_2(\mathbb F_q)$ distinguished for both $H_1$ and $K_1$, we attach an integer $r(\pi_1)$ as follows. For the principal series and cuspidal representations, respectively, we set
$r\bigl({\rm Ps}(\chi^j,\chi^{-j})\bigr)
 =j$ for $1\le j\le \frac{q-3}{2}$ and $r\bigl(\pi(\nu^{(q-1)(j+1)})\bigr)
  =j$ for $0\le j\le \frac{q-3}{2}$. For the remaining two representations, we set
$r(1)=0$ and $r({\rm St}\otimes\eta)=(q-1)/2$. Let $\pi_n\leftrightarrow\pi_1$ be as in Theorem \ref{thm1}. If $q=p^f$, let $\mathfrak o_f$ be the ring of integers of the unramified extension of $\mathbb Q_p$ of degree $f$.  

\begin{corollary}[Corollary to Theorem \ref{thm1}]\label{coro-thm1}
Fix an identification $\iota\colon\mathbb C\xrightarrow{\sim}\overline{\mathbb Q}_p$, and let $r=r_\iota(\pi_1)$ be the index attached to $\pi_1$ with respect to $\iota$. Then, for every $u\in\mathbb F_q$, the coefficient $\iota\bigl(c(\pi_n;H_n,K_{n,\alpha,u})\bigr)$ belongs to $\mathfrak o_f$, and
\[
\iota\bigl(c(\pi_n;H_n,K_{n,\alpha,u})\bigr)
\equiv \bigl[P_r(u\sqrt{\alpha})\bigr]^2\mod p,
\]
where $P_r(X)\in\mathbb Z[1/2][X]$ is the $r$-th Legendre polynomial.
\end{corollary}

We end the introduction by highlighting a conjecture originally made in \cite{ana25} which fits the theme of this paper via Corollary \ref{coro-thm1}. Its inspiration in \cite{ana25} came from the fact that the truth of this conjecture settles the missing case of the main result in \cite[\S 7.1]{bhko24} via the main result in \cite{ana25} in a neat way. When $q=p^f$ with $f>1$, \cite[Theorem 1.2]{ana25} gives a single $u\in\mathbb F_q^\times$ for which all the rank one correlation coefficients are $p$-adic units. By Theorem \ref{thm1} and Remark \ref{rmk-thm1}, the same $u$ works in every rank. The following conjecture would give the same conclusion in the remaining case $q=p$.

\begin{conjecture}\label{conj}
Let $Q_r(x)$ be such that 
\[
P_r(x) =
\begin{cases}
Q_r(x^2) &\text{if $r$ is even,} \\
xQ_r(x^2) &\text{if $r$ is odd,}  
\end{cases}
\]
where $P_r(X) \in \mathbb Z[\frac{1}{2}][X]$ is the $r$-th Legendre polynomial. Then
for every odd prime $p$, there exists a quadratic nonresidue $c$ modulo $p$ such that $Q_r(c)\not \equiv 0 \pmod p$ for every $0\le r\le(p-1)/2$. 
\end{conjecture}
 
\subsection*{Acknowledgements} 

The author would like to thank Dipendra Prasad for many helpful conversations and for his generous encouragement.

\subsection*{Use of AI} 
The author used GPT-6 Astra to search the literature, assist in developing and examining proofs, check calculations, review arguments for gaps, and draft portions of a first version of the manuscript. The author then verified all the arguments, calculations, and citations, and wrote the current manuscript. The final manuscript is run through GPT-6 Astra once again for minimal corrections. 

\section{Preliminaries}

In this section we fix the notations and collect the preliminary material. 

\subsection{Orthogonal groups and paramodular subgroups}\label{ortho-para}

We follow the exposition in \cite[\S\S 3,5]{gro16} and \cite[\S 2]{che22}.
Consider the lattice 
\[L_n= \mathbb Z e_1 \oplus \dots \oplus \mathbb Z e_n \oplus \mathbb Z v_0 \oplus \mathbb Z f_n \oplus \dots \oplus \mathbb Z f_1 \]
with the integral quadratic form
\[
Q_n\left(\sum_{i=1}^{n}x_i e_i+zv_0+
       \sum_{i=1}^{n}y_i f_i\right)
=z^2+\sum_{i=1}^{n}x_i y_i
\]
so that the associated symmetric bilinear form
\[
B_n(v,v^\prime)=Q_n(v+v^\prime)-Q_n(v)-Q_n(v^\prime)
\]
satisfies $B_n(v_0,v_0)=2$ and $B_n(e_i,f_i)= 1$ and all other pairings are zero.
Let $F$ be a finite or $p$-adic field and let 
\[
V_n = V_{n,F}=L_n \otimes_{\mathbb Z} F,
\]
and 
\[
G_n = SO(2n+1,F)
=\left\{
g \in {\rm SL}(V_{n,F}) \mid
Q_n(gv)=Q_n(v)~ \text{for all }v\in V_{n,F}
\right\}.
\]
To define the two even orthogonal subgroups of $G_n$ of interest in this paper, let $E/F$ be the quadratic extension when $F$ is finite 
and the unramified quadratic extension when $F$ is $p$-adic. 
Choose $\alpha\in F^\times\setminus F^{\times2}$, with $\alpha\in\mathfrak o_F^\times \setminus \mathfrak o_F^{\times 2}$ when $F$ is $p$-adic, such that $E=F(\sqrt\alpha)$.

We define
\[H_n = SO^+(2n,F) = {\rm Stab}_{G_n}(v_0),\] and 
\[K_{n,\alpha} = SO_\alpha^-(2n,F) = {\rm Stab}_{G_n} (e_n+\alpha f_n).\]

Now consider $n=1$. We identify our quadratic space $V_{1,F}$ with trace zero matrices under 
\[xe+zv_0+yf \longleftrightarrow 
\left(\begin{array}{cc} z & x \\ y & -z \end{array} \right), ~Q(X) = -\det X = z^2+xy.
\]
The action of ${\rm GL}_2(F)$ on this space by conjugation results in the isomorphism ${\rm PGL}_2(F) \simeq SO(3,F)$. Explicitly, this map is
\[
\left( \begin{array}{cc} a & b \\ c & d \end{array}\right) \mapsto \frac{1}{ad-bc} \left( 
\begin{array}{ccc}
a^2 & -2ab & -b^2 \\ -ac & ad+bc & bd \\ -c^2 & 2cd & d^2 
\end{array}
\right)
\]
and in particular the image of $[1,u;0,1]$ is  
\[
g^1_u = \left(
\begin{array}{ccc}
1 & -2u & -u^2 \\ 0 & 1 & u \\ 0 & 0 & 1
\end{array}
\right).
\]
Let $g_u = g^n_u = {\rm diag}(I_{n-1},g^1_u,I_{n-1})$. Define
\begin{equation}\label{k-def}
K_{n,\alpha,u} = g_u K_{n,\alpha} g_u^{-1}.
\end{equation}
Note that with the notation of \cite{ana25}, our rank-one torus $K_{1,\alpha,u}$ is $K_{\alpha^{-1},u}$.

We have the following elementary lemma whose proof we omit. 
\begin{lemma}\label{lemma-dc-f}
When $F=\mathbb F_q$, the elements $g_u$ form a set of representatives for $H_n\backslash G_n/K_{n,\alpha}$, as $u$ varies in $F$. In particular, $|H_n\backslash G_n/K_{n,\alpha}| = q$.
\end{lemma}

Now we define the relevant paramodular subgroups $P_{n,m}$ \cite[\S 2.2]{che22} (see also \cite[p. 684]{gro16}). Here, $F$ is a $p$-adic field and we denote by $\mathfrak o = \mathfrak o_F$ the ring of integers of $F$ and by $\mathfrak p = \varpi \mathfrak o$ its unique maximal ideal. Let $q= |\mathfrak o/\mathfrak p|$. Consider the lattice 
\[L_{n,m}= \mathfrak o e_1 \oplus \dots \oplus \mathfrak o e_n \oplus \mathfrak p^m v_0 \oplus \mathfrak p^m f_n \oplus \dots \oplus \mathfrak p^m f_1 \] and let
\[J_{n,m} = {\rm Stab}_{G_n}(L_{n,m}).\]
We define $P_{n,0} = J_{n,0} = SO(2n+1,\mathfrak o)$. For $m \geq 1$, let $D_m = {\rm diag}(I_n,\varpi^m,\varpi^m I_n)$. Since $L_{n,m} = D_m \mathfrak o^{2n+1}$, for $g \in G_n$, it follows that 
\[g \in J_{n,m} \iff g D_m\mathfrak o^{2n+1} = D_m\mathfrak o^{2n+1} \iff D_m^{-1}gD_m \in {\rm GL}_{2n+1}(\mathfrak o).\] 
Reducing $D_m^{-1}gD_m$ modulo $\mathfrak p$ and deleting its middle row and column, we get a morphism
\[\rho_m: J_{n,m} \rightarrow O^+(2n,\mathbb F_q).\]
Now the paramodular subgroup $P_{n,m}$ is the subgroup of $J_{n,m}$ of index $2$ with connected reduction; i.e., $P_{n,m}=\rho_m^{-1}(SO^+(2n,\mathbb F_q))$. As Gross remarks in \cite[p. 684]{gro16}, this definition is suggested by Brumer based on extensive computations in \cite{bk14}. 
\begin{remark}\label{rmk-para-conj}
With $b_{n,r}={\rm diag}(\varpi^r I_n,1,\varpi^{-r}I_n)$ and
\[
C_n(r)=\ker \left(G_n(\mathfrak o)\longrightarrow
G_n(\mathfrak o/\mathfrak p^r)\right),
\]
we have, by \cite[Proposition 7.1.3]{tsa13} (or \cite[Corollary 3.3]{che26}), 
\[b_{n,r}P_{n,2r}b_{n,r}^{-1}=H_n(\mathfrak o)C_n(r).\]
For $j\ge1$ and $m\ge0$, we define
\[
\mathcal P_{j,m}
=
b_{j,\lfloor m/2\rfloor}P_{j,m}
b_{j,\lfloor m/2\rfloor}^{-1}.
\]
This agrees with Cheng’s $K^0_{j,m}$ \cite[(3.8)]{che26}.
\end{remark}

\begin{remark}
Under the explicit identification between ${\rm PGL}_2(F)$ and $SO(3,F)$ stated above, we have $P_{1,m} = \Gamma_0(\mathfrak p^m)$, the Atkin-Lehner congruence subgroup of ${\rm GL}_2(F)$.
Similarly, for $n=2$, starting with $D=F^4$ with symplectic basis $d_1,d_2,d_3,d_4$, 
$\langle d_1,d_4\rangle=\langle d_2,d_3\rangle=1$, we identify $V_{2,F}$ with the quadratic space 
\[
\ker \left(\bigwedge\nolimits^2D\longrightarrow F,\quad
u\wedge v\longmapsto\langle u,v\rangle\right)
\]
by
\[
(e_1,e_2,v_0,f_2,f_1)
=(d_{12},d_{13},d_{14}-d_{23},d_{24},-d_{34}),
\]
where $d_{ij}=d_i\wedge d_j$ and the quadratic form $Q$ is defined by
\[
z\wedge z=-2Q(z)\,d_1\wedge d_2\wedge d_3\wedge d_4,
\]
for $z \in \ker\!\left(\bigwedge\nolimits^2D\longrightarrow F\right)$. 
Note that
\[
Q(x_1e_1+x_2e_2+zv_0+y_2f_2+y_1f_1)
=z^2+x_1y_1+x_2y_2.
\]
We get ${\rm PGSp}_4(F) \simeq SO(5,F)$ obtained via 
restricting $\mu(g)^{-1}\bigwedge^2g$ to this kernel, where $\mu(g)$ is the symplectic similitude factor. Under this explicit identification, we have $P_{2,m}=K(\mathfrak p^m)$, the Roberts-Schmidt paramodular subgroup of ${\rm GSp}_4(F)$ \cite[p. 6]{rs07}. See \cite[Chapter 6]{tsa13} for further details. Note that $\Gamma_0(\mathfrak p^m)$ and $K(\mathfrak p^m)$ here denote the images of the standard congruence subgroups in the corresponding projective groups.
\end{remark}

\subsection{Theta correspondence over finite and $p$-adic fields}\label{theta-finite}

We broadly follow the exposition in \cite{kud96,pra98}. Let $F$ be a finite or $p$-adic field. Let $\psi$ be a nontrivial unitary additive character of $F$. Let $W=F^2$ with symplectic basis $w^+, w^-, \langle w^+, w^- \rangle =1$ and let $(X,Q_X)$ be a nondegenerate quadratic space over $F$. Let $\omega_{X,\psi}$ denote the corresponding Weil representation which is a representation of ${\rm Sp}(W \otimes X)$ (resp. ${\rm Mp}(W\otimes X)$) when $F$ is finite (resp. $p$-adic). 
 
Now we recall the Schr\"odinger model for $\omega_{X,\psi}$. Let $d=\dim X$ and 
\[
B_X(x,y)=Q_X(x+y)-Q_X(x)-Q_X(y).
\]
We equip $W\otimes X$ with the symplectic form
\[
\langle w\otimes x,w'\otimes y\rangle
=\langle w,w'\rangle B_X(x,y).
\]
The Lagrangian decomposition
\[
W\otimes X=(w^+\otimes X)\oplus(w^-\otimes X)
\]
realizes the Weil representation $\omega_{X,\psi}$ of
${\rm Sp}(W\otimes X)$ (or its metaplectic cover in the $p$-adic case) on the inner product space of complex valued functions 
$S(w^-\otimes X)\simeq\mathcal S(X)$ where 
\[
\mathcal S(X)=
\begin{cases}
\mathbb C[X], &\text{$F$ finite},\\
C_c^\infty(X), &\text{$F$ $p$-adic},
\end{cases}
\]
and 
\[
\langle f_1,f_2\rangle=
\begin{cases}
\displaystyle\sum_{x\in X}\overline{f_1(x)}f_2(x),
&\text{$F$ finite},\\
\displaystyle\int_X\overline{f_1(x)}f_2(x)\,d_\psi x,
&\text{$F$ $p$-adic},
\end{cases}
\]
where $d_\psi x$ is the self-dual measure for the pairing
$(x,y)\mapsto\psi(B_X(x,y))$.

We write $\omega_{X,\psi}(s,h)$ for its restriction through
\[
{\rm SL}_2(F)\times SO(X)
\longrightarrow {\rm Sp}(W\otimes X), ~
(s,h) \mapsto s\otimes h
\]
or in the $p$-adic case its standard metaplectic lift 
\[{\rm Mp}_2(F) \times SO(X) \longrightarrow {\rm Mp}(W \otimes X).\]
The action of ${\rm Mp}_2(F)$ factors through ${\rm SL}_2(F)$ precisely when $d$ is even.

Define
\[
\widehat f(x)=\int_X f(y)\psi(-B_X(x,y))\,d_\psi y.
\]
If $F=\mathbb F_q$, this means
\[
\widehat f(x)
=q^{-d/2}\sum_{y\in X}\psi(-B_X(x,y))f(y).
\]
Let $\gamma_{X,\psi}$ be the Weil index of $Q_X$. When $F=\mathbb F_q$, 
\[
\gamma_{X,\psi}
=q^{-d/2}\sum_{x\in X}\psi(Q_X(x)).
\]
For
$b\in F$ and $a\in F^\times$,
let
\[
u(b)=\begin{pmatrix}1&b\\0&1\end{pmatrix},~
m(a)=\begin{pmatrix}a&0\\0&a^{-1}\end{pmatrix},~
j=\begin{pmatrix}0&-1\\1&0\end{pmatrix},
\]
and we denote by the same symbols their compatible standard lifts in the metaplectic case. The actions of these elements on the Schr\"odinger model are given by the following formulas (cf.~\cite[pp.~64--65]{ger77},
\cite[\S 1.1 and \S 4]{pra98}, and
\cite[\S II.4]{kud96}):
\begin{align*}
[\omega_{X,\psi}(u(b),1)f](x)
&=\psi(bQ_X(x))f(x),\\
[\omega_{X,\psi}(m(a),1)f](x)
&=\frac{\gamma_{X,\psi}}{\gamma_{X,\psi_a}}
  |a|^{d/2}f(ax),\\
[\omega_{X,\psi}(j,1)f](x)
&=\gamma_{X,\psi}^{-1}\widehat f(x),\\
[\omega_{X,\psi}(1,h)f](x)
&=f(h^{-1}x).
\end{align*}
Here, $|a|$ is the normalized absolute value in the
$p$-adic case and is defined to be $1$ in the
finite-field case. The character $\psi_a$ is defined by $\psi_a(y)=\psi(ay)$.
If $F$ is finite and $\varepsilon$ is its
quadratic character, then
\[
\frac{\gamma_{X,\psi}}{\gamma_{X,\psi_a}}
=\varepsilon(a)^d.
\]

This construction is compatible with orthogonal direct sums.
If $X=X_1\perp X_2$, the natural identification
\[
\mathcal S(X)=\mathcal S(X_1)\otimes\mathcal S(X_2)
\]
satisfies
\[
\omega_{X,\psi}(s,h_1\oplus h_2)
=
\omega_{X_1,\psi}(s,h_1)\otimes
\omega_{X_2,\psi}(s,h_2).
\]
This is the tensor product property of the Weil representation
\cite[\S 1, Remark 1]{pra98} (see also \cite[Corollary 2.5]{ger77}).

Now, let $N=\{u(b) \mid b\in F\}$, identified with its canonical lift
in the metaplectic case, and $\psi_t(u(b))=\psi(bt)$. Define the twisted coinvariant space by
\[
\mathcal S(X)_{N,\psi_t}
=
\mathcal S(X)\Big/
{\rm span}\{
\omega_{X,\psi}(u(b),1)f-\psi(bt)f \mid 
b\in F, f \in \mathcal S(X)\}.
\]
From the formula for $\omega_{X,\psi}(u(b),1)$, we get the following lemma (see the proof of \cite[ Proposition 6.1]{gw21}).

\begin{lemma}\label{lemma-norm}
For every $t\in F$, restriction to
$X_t=\{x\in X \mid Q_X(x)=t\}$ induces an
$\operatorname{SO}(X)$-equivariant isomorphism
\[
\mathcal S(X)_{N,\psi_t}\simeq\mathcal S(X_t).
\]
Here, $\mathcal S(X_t)$ denotes all functions on $X_t$ when
$F$ is finite, and the locally constant compactly supported
functions on $X_t$ when $F$ is $p$-adic.
In particular, when $F$ is finite, extension by zero identifies
\[
\mathbb C[X]^{N,\psi_t}=\mathbb C[X_t].
\]
\end{lemma}

\begin{remark}\label{rmk-norm}
In the notation for $V_n$, $H_n$, and
$K_{n,\alpha}$ as in the finite field case of \S 2.1,
\[
\mathbb C[V_n]^{N,\psi_1}
\simeq {\rm Ind}_{H_n}^{G_n} 1 ~\&~
\mathbb C[V_n]^{N,\psi_\alpha}
\simeq {\rm Ind}_{K_{n,\alpha}}^{G_n} 1.
\]
\end{remark}

We apply the tensor product property in our situation as follows.
In the notation of \S\ref{ortho-para}, let
\begin{equation}\label{directsum}
T_n=\langle e_n,v_0,f_n\rangle,~
U_{n-1}=\bigoplus_{i=1}^{n-1}\langle e_i,f_i\rangle,
\end{equation}
so that $V_n=U_{n-1}\perp T_n$, with $T_n\simeq V_1$.
Writing
\[
W_1=W\otimes U_{n-1},~
W_2=W\otimes T_n\simeq W\otimes V_1,
\]
we obtain the orthogonal symplectic decomposition
\[
W\otimes V_n=W_1\oplus W_2.
\]
The tensor product property of the Weil representation,
applied to the symplectic groups in the finite case
and their metaplectic covers in the $p$-adic case,
gives the natural identification
\begin{equation}\label{eq-tensor}
\mathcal S(V_n)
\simeq\mathcal S(U_{n-1})\otimes\mathcal S(V_1).
\end{equation}
Under this identification,
\begin{equation}\label{tensor}
\omega_{n,\psi}(s,i_n(g))
=
\omega_{U_{n-1},\psi}(s,1)
\otimes\omega_{1,\psi}(s,g),~
s\in M_F, ~g\in SO(V_1),
\end{equation}
where 
\[
M_F=
\begin{cases}
{\rm SL}_2(F),&F\text{ finite},\\
{\rm Mp}_2(F),&F\text{ $p$-adic},
\end{cases}
\]
and
$i_n:SO(V_1)\hookrightarrow SO(V_n)$ acts trivially
on $U_{n-1}$.
Since $\dim U_{n-1}$ is even, its Weil representation
factors through ${\rm SL}_2(F)$ also in the $p$-adic case.

To define the big theta lift, we restrict the Weil representation
$\omega_{n,\psi}$ to $M_F\times SO(V_n)$.
For an irreducible representation $\sigma$ of $M_F$
(smooth and genuine in the $p$-adic case), the maximal
$\sigma$-isotypic quotient of $\omega_{n,\psi}$ has the form
\begin{equation}\label{eq-theta}
(\omega_{n,\psi})_\sigma
\simeq \sigma\otimes\Theta_{n,\psi}(\sigma).
\end{equation}
This defines the big theta lift $\Theta_{n,\psi}(\sigma)$,
a possibly zero representation of $G_n = SO(V_n)=SO(2n+1,F)$
(cf.~\cite[\S 4]{pra98}).
Equivalently,
\begin{equation}\label{eq-can}
\Theta_{n,\psi}(\sigma)
=\bigl(\mathcal S(V_n)\otimes\sigma^\vee\bigr)_{M_F},
\end{equation}
where the subscript denotes coinvariants for the $M_F$-action
$s \mapsto \omega_{n,\psi}(s,1)\otimes\sigma^\vee(s)$.
When $F$ is finite, by semisimplicity, we have the decomposition
\[
\omega_{n,\psi}|_{{\rm SL}_2(F)\times SO(V_n)}
\simeq
\bigoplus_{\sigma\in{\rm Irr}({\rm SL}_2(F))}
\sigma\otimes\Theta_{n,\psi}(\sigma).
\]

For $p$-adic $F$, by Howe duality, $\Theta_{n,\psi}(\sigma)$
has finite length and, when nonzero, has a unique irreducible
quotient, which we denote by $\theta_{n,\psi}(\sigma)$. This is the small theta lift. Moreover, for irreducible genuine smooth representations
$\sigma$ and $\sigma^\prime$ of ${\rm Mp}_2(F)$,
\[
\theta_{n,\psi}(\sigma)\simeq\theta_{n,\psi}(\sigma^\prime)\ne0
\Longrightarrow
\sigma\simeq\sigma^\prime.
\]
See \cite[\S 4, Theorem 6]{pra98} and
\cite[\S II.2, Howe Duality Principle]{kud96}. If $\sigma$ is supercuspidal, Kudla's theorem asserts that
$\Theta_{n,\psi}(\sigma)$ is either zero or irreducible \cite[\S III.6, Theorem 6.1(i)]{kud96},
and therefore
\[
\Theta_{n,\psi}(\sigma)\simeq\theta_{n,\psi}(\sigma)
\]
in this case. Let 
$\theta_{n,\psi}^{\leftarrow}(\pi)$ denote the small theta lift
from $G_n$ to ${\rm Mp}_2(F)$. Let
\begin{align*}
{\rm Howe}_\psi({\rm Mp}_2(F);G_n)
&=\{\sigma\in {\rm Irr}_{{\rm gen}}({\rm Mp}_2(F))
       \mid \theta_{n,\psi}(\sigma)\neq 0\},\\
{\rm Howe}_\psi(G_n;{\rm Mp}_2(F))
&=\{\pi\in {\rm Irr}(G_n)
       \mid \theta_{n,\psi}^{\leftarrow}(\pi)\neq 0\},
\end{align*}
where ${\rm Irr}_{{\rm gen}}({\rm Mp}_2(F))$ denotes the
isomorphism classes of irreducible genuine smooth
representations of ${\rm Mp}_2(F)$. Equivalently, these are the sets
of representations satisfying
\[
{\rm Hom}_{{\rm Mp}_2(F)}(\mathcal S(V_n),\sigma)\neq 0,
~
{\rm Hom}_{G_n}(\mathcal S(V_n),\pi)\neq 0,
\]
respectively. The Howe duality principle then asserts that
\begin{equation}\label{eq-howe}
{\rm Howe}_\psi({\rm Mp}_2(F);G_n)
\xrightarrow{\;\sim\;}
{\rm Howe}_\psi(G_n;{\rm Mp}_2(F)),
~
\sigma\longmapsto\theta_{n,\psi}(\sigma),
\end{equation}
is a bijection, where the inverse is given by
$\pi\longmapsto\theta_{n,\psi}^{\leftarrow}(\pi)$;
see \cite[\S II.2, p.~33]{kud96}.

For an irreducible admissible genuine representation $\sigma$ of ${\rm Mp}_2(F)$, we know that if $\pi_j=\theta_{j,\psi}(\sigma)\neq0$ for some $j\geq1$, then $\pi_n=\theta_{n,\psi}(\sigma)\neq 0$ for every $n\geq j$. This is the so called persistence principle \cite[Chapter III, \S 4, Proposition 4.1(i)]{kud96}. Moreover, $\pi_n$ is a subquotient of a representation obtained from $\pi_j$ by successive normalized parabolic induction with unramified characters of ${\rm GL}_1(F)$ \cite[Chapter III, \S 4, Proposition 4.1(iii)]{kud96}. If, in addition, $\pi_j$ is also tempered, then this can be strengthened further as: for every $n>j$, $\pi_n$ is the unique irreducible quotient of the standard module
\[
|\cdot|^{n-\frac32}\times\cdots\times
|\cdot|^{j-\frac12}\rtimes\pi_j,
\]
where also the induction is normalized \cite[Proposition 3.2]{gt16}.

The most important case for us in this paper is the base case $n=1$, where our theta correspondence is the Waldspurger correspondence \cite{wal80,wal91}. We state only what we need and refer to \cite[\S 5]{gs12} for a concise summary of all the other details. 
Here, $D$ is the quaternion division algebra over $F$, and $PD^\times=D^\times/F^\times$. By its conjugation action on the trace zero subspace $D^0$, we identify $PD^\times$ with $SO(D^0,-{\rm Nrd})$, the nonsplit form of $SO(3)$. For the second part, see \cite[Theorem 8.1(ii)]{gs12}.

\begin{proposition}\label{prop-waldspurger}
For every irreducible
admissible representation $\tau$ of ${\rm PGL}_2(F)$ or $PD^\times$,
its small theta lift $\theta_{1,\psi}^{\leftarrow}(\tau)$ to ${\rm Mp}_2(F)$ is
nonzero (and irreducible). Moreover, the assignment
$\tau\mapsto \theta_{1,\psi}^{\leftarrow}(\tau)$ defines a bijection
\[
{\rm Irr}({\rm PGL}_2(F))
\sqcup
{\rm Irr}(PD^\times)
\xrightarrow{\;\sim\;}
{\rm Irr}_{\rm gen}({\rm Mp}_2(F)).
\]
Also, if $\tau$ is tempered, its big theta lift is irreducible and
\[
\Theta^{\leftarrow}_{1,\psi}(\tau)
\simeq
\theta^{\leftarrow}_{1,\psi}(\tau).
\]
\end{proposition}

We end this section by stating the following result for finite fields for which our main reference is the work of Kriz (see \cite[\S 6]{kri25} for the first part and \cite[Theorem 1 (2)]{kri25} or \cite[Theorem 2 (2)]{kri24} for the second part) but for $n\ge3$, both assertions also follow from \cite[Proposition 5.7 ~\&~ Table I]{lost10}.
\begin{proposition}\label{kriz}
Let $\sigma$ be an irreducible representation of ${\rm SL}_2(F)$ which extends to an irreducible representation of ${\rm GL}_2(F)$. Then,
\begin{enumerate}
\item $\dim \Theta_{n,\psi}(\sigma) = \frac{q^{2n}-1}{q^2-1}\dim \sigma$,
\item $\Theta_{n,\psi}(\sigma)$ is an irreducible representation of $SO(2n+1,F)$.
\end{enumerate}
\end{proposition}

\begin{remark}
In Proposition \ref{kriz} (2), to appeal to \cite{kri24,kri25} we need the stable range assumption which is $n \geq 2$. For nontrivial $\sigma$, we give proofs of (1) and (2), including the case $n=1$, in §\ref{kriz1} and \S \ref{kriz2}. When $n=1$ and $\sigma=1$, the theta lift is one dimensional and hence irreducible.
\end{remark}

\section{Proof of Theorem \ref{thm1}}\label{proof-thm1}

In this section we prove Theorem \ref{thm1} in \S \ref{proofthm1}. The basic idea behind the proof is simple. When there is an ambient representation which may be large but then concrete, say a function space, due to the largeness though we lose multiplicity one for invariant vectors but due to concreteness natural choices for invariant vectors may be available. It is typically easy to do calculations in such a concrete situation and these calculations are then transferred to irreducible constituents by projection and Schur orthogonality, after establishing nonvanishing of the projected vectors. A similar viewpoint is made use of in \cite[Corollary 1.2 and Remark 9]{am20}, where a natural invariant function in the Gelfand-Graev representation projects nontrivially to every generic distinguished constituent. 

The natural vectors for us here are the delta functions at $v_0$ and $e_n+\alpha f_n$ in the
Schr\"odinger model of the Weil representation $\omega_{n,\psi}$ of ${\rm Sp}(W \otimes V_n)$. By
making crucial use of the tensor identity (\ref{tensor}) in \S \ref{theta-finite}, we reduce their matrix
coefficients to rank one. Multiplicity one for the relevant Whittaker spaces of $\sigma$ of 
${\rm SL}_2(\mathbb F_q)$ is then used in the $\sigma$-isotypic 
component $\sigma\otimes\Theta_{n,\psi}(\sigma)$ to express the projected delta functions as pure tensors $e_\pm\otimes w_{n,\pm}$. All these are carried out in \S \ref{two lemmas}.

The matrix coefficient between these vectors therefore factors into a matrix coefficient of $\sigma$ and one of $\Theta_{n,\psi}(\sigma)$. We make use of irreducibility and dimension formula for the theta lifts, due to Kriz \cite{kri24,kri25}, stated as Proposition \ref{kriz}, to identify the distinguished representations and for the normalization in Theorem \ref{thm1}. We give independent proofs of these results in \S \ref{kriz1} and \S \ref{kriz2}. We can then deduce Theorem \ref{thm1} through Schur orthogonality which we do in \S \ref{proofthm1}.

\subsection{Two lemmas}\label{two lemmas}
 
We take $X = V_n$ in the notation of \S \ref{theta-finite}. Recall that for $g\in G_n$ the Weil representation action is given by
\[
[\omega_{n,\psi}(1,g)f](x)=f(g^{-1}x).
\]
Thus, for the characteristic function $\delta_v^n$ of $\{v\}$,
\[
\omega_{n,\psi}(1,g)\delta_v^n = \delta_{gv}^n.
\]
In particular, the unit vectors $\delta_{v^+}^n$ and $\delta_{v^-}^n$, where $v^+=v_0$ and $v^- = e_n+\alpha f_n$, are fixed respectively by $H_n$ and $K_{n,\alpha}$.

Under the identification (cf. (\ref{directsum})),
\begin{equation}\label{eq-pure-tensor}
\mathbb C[V_n]
\simeq
\mathbb C[U_{n-1}]\otimes\mathbb C[T_n],~
(\phi\otimes \xi)(x+y)=\phi(x)\xi(y),
\end{equation}
both these characteristic functions factor as a tensor product of characteristic functions
\[
\delta_{v^\pm}^{\,n}
=
\delta_0^{U_{n-1}}\otimes\delta_{v^\pm}^{\,1},
\]
where for the second factor we have made use of the natural identification $T_n\simeq V_1$.

Now, for $s = [a,b;c,d] \in {\rm SL}_2(\mathbb F_q)$, let 
\[
\omega_{U,\psi}(s)=\omega_{U_{n-1},\psi}(s,1),
\]
where $(s,1)$ maps to
\[
s \otimes I_{U_{n-1}}
=
\begin{pmatrix}
aI_{2n-2}&bI_{2n-2}\\
cI_{2n-2}&dI_{2n-2}
\end{pmatrix}
\in {\rm Sp}(W_1).
\]
Also, for $g_1\in G_1$, let
\[
i_n(g_1) = I_{U_{n-1}} \oplus g_1
=
\begin{pmatrix}
I_{n-1}&0&0\\
0&g_1&0\\
0&0&I_{n-1}
\end{pmatrix}
\in SO(V_n),
\]
using $T_n\simeq V_1$.
The tensor product identity (cf. (\ref{tensor})) gives
\begin{equation}\label{eqn-kappa}
\left\langle
\delta_{v^+}^n,
\omega_{n,\psi}(s,i_n(g_1))\delta_{v^-}^n
\right\rangle 
=
\kappa_{n-1}(s)
\left\langle
\delta_{v^+}^1,
\omega_{1,\psi}(s,g_1)\delta_{v^-}^1
\right\rangle,
\end{equation}
where
\[
\kappa_{n-1}(s)
=
\left\langle
\delta_0^{U_{n-1}},
\omega_{U,\psi}(s)\delta_0^{U_{n-1}}
\right\rangle.
\]
To compute this scalar, take $X = U_{n-1}$ and apply the Schr\"odinger formulas of \S \ref{theta-finite}. Observe that $u(b)$ and $m(a)$ fix
$\delta_0 = \delta_0^{U_{n-1}}$, while
\[
\omega_{U,\psi}(j)\delta_0
=
q^{-(n-1)}.
\]
For $s = [a,b;c,d] \in {\rm SL}_2(\mathbb F_q)$, 
using the decomposition
\[
s=
\begin{cases}
u(ab)m(a),&c=0,\\
u(a/c)m(c^{-1})j\,u(d/c),&c\ne0,
\end{cases}
\]
evaluating $\omega_{U,\psi}(s)\delta_0$ at $0$ gives
\[
\kappa_{n-1}(s) =
\begin{cases}
1,&c=0,\\
q^{-(n-1)},&c\ne0.
\end{cases}
\]
Note that when $c=0$, the function $\omega_{1,\psi}(s,g_1)\delta_{v^-}^1$
is supported at $a^{-1}g_1v^-$, and 
$Q_1(a^{-1}g_1v^-) = \alpha/a^2$. Since $\alpha$ is a nonsquare, this can never be $Q_1(v^+) = 1$. Hence, from (\ref{eqn-kappa}) we conclude the following lemma. 
\begin{lemma}\label{lemma-weil} 
For every $s \in {\rm SL}_2(\mathbb F_q)$ and $g \in G_1$, we have
\[
\left\langle
\delta_{v^+}^n,
\omega_{n,\psi}(s,i_n(g))\delta_{v^-}^n
\right\rangle 
=
q^{-(n-1)}
\left\langle
\delta_{v^+}^1,
\omega_{1,\psi}(s,g)\delta_{v^-}^1
\right\rangle.
\]
\end{lemma}

Now from Lemma \ref{lemma-weil} we want to conclude a similar relation between the matrix coefficients 
$\bigl(v_{H_j},\Theta_{j,\psi}(\sigma)(i_j(g_1))v_{K_{j,\alpha}}\bigr)$
for $j=1$ and $j=n$. Here, $v_\square$ denotes a fixed vector for the indexed subgroup $\square$ in the representation $\Theta_{j,\psi}(\sigma)$ (see (\ref{eq-theta}) for the definition) where $\sigma$ is any nontrivial irreducible representation of ${\rm SL}_2(\mathbb F_q)$ that extends to ${\rm GL}_2(\mathbb F_q)$ as in Proposition \ref{kriz}; thus there are $(q-1)$ many such representations. It is well known that for each such $\sigma$, we have $\dim \sigma^{N,\psi}=1$ for any nontrivial additive character $\psi$ of $\mathbb F_q$.

As a first step towards this, we consider the orthogonal projection $P_{n,\sigma}$ from the Weil representation $\omega_{n,\psi}$ onto its $\sigma$-isotypic component for the ${\rm SL}_2(\mathbb F_q)$-action:
\[
P_{n,\sigma}\mathbb C[V_n]
\simeq \sigma\otimes\Theta_{n,\psi}(\sigma) 
\]
given by the formula
\begin{equation}\label{proj}
P_{n,\sigma}
=\frac{\dim\sigma}{|{\rm SL}_2(\mathbb F_q)|}
\sum_{s\in {\rm SL}_2(\mathbb F_q)}
\overline{\chi_\sigma(s)} \omega_{n,\psi}(s,1),
\end{equation}
and we show that $P_{n,\sigma}(\delta^n_{v^\pm})$ are nonzero. 
\begin{lemma}\label{lemma-nonzero}
Let $v \in V_n$ with $t = Q_n(v) \neq 0$. Then,
\[
\|P_{n,\sigma}\delta_v^n\|^2
=
\frac{\dim\sigma^{N,\psi_t}
      \dim\Theta_{n,\psi}(\sigma)}
     {[G_n: {\rm Stab}_{G_n}(v)]}
=  \frac{\dim\Theta_{n,\psi}(\sigma)}
     {q^n(q^n+\varepsilon(t))}.
\]
Consequently, $P_{n,\sigma}\delta_{v^\pm}^n\neq0$ and choosing unit generators $e_\pm$ in $(N,\psi_t)$-eigenspaces, where $t=1,\alpha$, we get 
\[
P_{n,\sigma}\delta_{v^\pm}^n= e_\pm\otimes w_\pm
\]
with
\[
0 \neq w_+ \in\Theta_{n,\psi}(\sigma)^{H_n},~
0 \neq w_-\in\Theta_{n,\psi}(\sigma)^{K_{n,\alpha}}.
\]
\end{lemma}
\begin{proof}
The restriction of $P_{n,\sigma}$ to $\mathbb C[V_n]^{N,\psi_t}$ is a $G_n$-equivariant orthogonal projection with image
\[
P_{n,\sigma} (\mathbb C[V_n]^{N,\psi_t})
\simeq
\sigma^{N,\psi_t}\otimes\Theta_{n,\psi}(\sigma).
\]
By Remark \ref{rmk-norm}, the lemma follows from the elementary observation: 
for any finite group $G$, subgroup $H$, and $G$-equivariant orthogonal projection $P$ on $\mathbb C[G/H]$,
\[
{\rm rank~} P
= {\rm Tr~} P
=\sum_{x\in G/H} \langle \delta_x,P\delta_x \rangle
=\sum_{x\in G/H}\|P\delta_x\|^2
=[G:H] \|P\delta_H\|^2,
\]
where the last equality follows from transitivity and equivariance. Invariance for $w_\pm$ follows from respective invariance for $\delta^n_{v^\pm}$ (since $P_{n,\sigma}$ is equivariant for $G_n$ and hence for its subgroups).

To see the cardinality of the index in the denominator, fix $v\in V_n$ with $Q_n(v)=t\ne0$. We know that $[G_n:\operatorname{Stab}_{G_n}(v)]$ is the number of vectors in the $G_n$-orbit of $v$. This orbit is the entire fiber $Q_n^{-1}(t)$. Indeed, for any $w$ with $Q_n(w)=t$, the map $v\mapsto w$ is an isometry of their spans and Witt’s theorem extends it to an isometry of $V_n$. 
If the extension has determinant $-1$, choose $u\in w^\perp$ with $Q_n(u)\ne0$ and compose the extension with the reflection in $u$. Thus, by \cite[Theorem 6.27]{ln97}, 
\[
[G_n:\operatorname{Stab}_{G_n}(v)]
=|\{v \mid Q_n(v)=t\}|
=q^n\bigl(q^n+\varepsilon(t)\bigr).
\] 
\end{proof}

\subsection{Dimension formula}\label{kriz1}

We can now give an elementary proof of Proposition \ref{kriz} (1) due to \cite{lost10,kri25}. Note that
\[\|P_{n,\sigma}\delta^n_v\|^2 = \langle P_{n,\sigma}\delta^n_v, P_{n,\sigma}\delta^n_v \rangle = \langle \delta^n_v, P_{n,\sigma}\delta^n_v \rangle = \frac{\dim\sigma}{|{\rm SL}_2(\mathbb F_q)|}
\sum_{s\in {\rm SL}_2(\mathbb F_q)}
\overline{\chi_\sigma(s)} \langle \delta^n_v, \omega_{n,\psi}(s,1)\delta^n_v \rangle. \]
The presence of the matrix coefficient
\[k_{n,v}(s) =  \langle \delta^n_v, \omega_{n,\psi}(s,1)\delta^n_v \rangle\]
is indicative of an $n$ to $1$ reduction as in Lemma \ref{lemma-weil} which we now do. 
Take $v=v^\pm$, and identify it with the corresponding vector of $V_1$.
Writing 
\[\delta^n_v = \delta_0^{U_{n-1}} \otimes \delta^1_v\]
as in the proof of Lemma \ref{lemma-weil} and proceeding as in the next steps there, we get
\begin{align*}
k_{n,v}(s) &=
\begin{cases}
k_{1,v}(s),& c=0,\\
q^{-(n-1)}k_{1,v}(s),&c \neq 0.
\end{cases} \\
&= [q^{-(n-1)}+(1-q^{-(n-1)})1_{B(\mathbb F_q)}] k_{1,v}(s),
\end{align*}
where $B$ is the upper triangular Borel subgroup of ${\rm SL}_2$.
Unlike in the proof of Lemma \ref{lemma-weil}, we may have nonzero value even at $c=0$. In any case, we have
\[\|P_{n,\sigma}\delta^n_v\|^2 = q^{-(n-1)} \|P_{1,\sigma}\delta^1_v\|^2 + (1-q^{-(n-1)})  \frac{\dim\sigma}{|{\rm SL}_2(\mathbb F_q)|}
\sum_{s\in B(\mathbb F_q)}
\overline{\chi_\sigma(s)} k_{1,v}(s). \]
Note that, for $s=[a,b;0,a^{-1}]$ and $t=Q_1(v) \neq 0$,
\[
k_{1,v}(s)
=[\omega_{1,\psi}(s,1)\delta_v](v)
=\varepsilon(a)\psi(abt)\,\delta_v(av) =
\begin{cases}
\psi(bt),&a=1,\\
0,&a \neq 1.
\end{cases}
\]
Thus, the Borel sum above is 
\[q \cdot \langle \chi_\sigma, \psi_t \rangle_N = q\] as our representations $\sigma$ afford Whittaker models for all nontrivial additive characters.
By Lemma \ref{lemma-nonzero}, 
\[
 \frac{\dim\Theta_{n,\psi}(\sigma)} {q^n(q^n+\varepsilon(t))}
 = q^{-(n-1)}  \frac{\dim\Theta_{1,\psi}(\sigma)}
     {q(q+\varepsilon(t))} + (1-q^{-(n-1)}) \frac{\dim \sigma}{q^2-1}.
\]
With $D_j=\dim\Theta_{j,\psi}(\sigma)$ and $d_\sigma=\dim\sigma$, write this as,
\[\frac{D_n}{q^n \pm 1} = \frac{D_1}{q \pm 1} + \frac{(q^n-q)d_\sigma}{q^2-1}.\]
Subtracting one from the other, it follows that 
\[D_n = \frac{q^{2n}-1}{q^2-1}D_1.\] Substituting this back,
\[(q^n-q)(D_1-d_\sigma)=0\] and thus $D_1=d_\sigma$ as well.

\begin{remark}\label{rmk-normsquared}
In the notations of Lemma \ref{lemma-nonzero}, 
\[
\|w_\pm\|^2
=\left\|P_{n,\sigma}\delta_{v^\pm}^n\right\|^2
=\frac{\dim\Theta_{n,\psi}(\sigma)}
       {q^n(q^n\pm1)}.
\]
Thus, making use of the dimension formula (see Proposition \ref{kriz} (1) and \S \ref{kriz1})
\[
\dim\Theta_{n,\psi}(\sigma)
=\dim\sigma\,\frac{q^{2n}-1}{q^2-1},
\]
we obtain
\[
\|w_+\|^2=\frac{\dim\sigma}{q^2-1}(1-q^{-n}),~
\|w_-\|^2=\frac{\dim\sigma}{q^2-1}(1+q^{-n}).
\]
\end{remark}

\subsection{Irreducibility and multiplicity one}\label{kriz2}

We now give an elementary proof of Proposition \ref{kriz} (2) due to \cite{lost10,kri24}. Multiplicity one for representations distinguished by both $H_n$ and $K_{n,\alpha}$ also follows from its proof.

Let $\mathcal R$ be our family of $q-1$ nontrivial 
irreducible representations of ${\rm SL}_2(\mathbb F_q)$ as in Proposition \ref{kriz}.

\begin{lemma}\label{lemma-kriz2}
For every $n\geq 1$, the representations
$\Theta_{n,\psi}(\sigma)$, $\sigma\in\mathcal R$, are nontrivial,
irreducible, and pairwise nonisomorphic. Moreover,
\[
\dim\Theta_{n,\psi}(\sigma)^{H_n}
=
\dim\Theta_{n,\psi}(\sigma)^{K_{n,\alpha}}
=1.
\]
Together with $1_{G_n}$, they exhaust the irreducible representations
of $G_n$ distinguished by both $H_n$ and $K_{n,\alpha}$.
\end{lemma}

\begin{proof}
Let (cf. Lemma \ref{lemma-norm})
\[
\mathcal M_+
=\mathbb C[\{x\in V_n:Q_n(x)=1\}]
\simeq \operatorname{Ind}_{H_n}^{G_n}1,
~
\mathcal M_-
=\mathbb C[\{x\in V_n:Q_n(x)=\alpha\}]
\simeq \operatorname{Ind}_{K_{n,\alpha}}^{G_n}1.
\]
As 
\[
\mathbb C[\{x\in V_n:Q_n(x)=t\}]
\simeq
\bigoplus_{\tau \in {\rm Irr}({\rm SL}_2(\mathbb F_q))}
\tau^{N,\psi_t}\otimes\Theta_{n,\psi}(\tau),
\]
it follows from multiplicity one for Whittaker models for ${\rm SL}_2(\mathbb F_q)$ and Lemma \ref{lemma-norm} that
\[
\bigoplus_{\sigma\in\mathcal R}\Theta_{n,\psi}(\sigma)
\]
is a direct summand of each $\mathcal M_\pm$.

Observe that none of these summands contains $1_{G_n}$. Indeed, for the 
${\rm SL}_2(\mathbb F_q)$-action,
\[
\mathbb C[V_n]^{G_n}
\simeq
{\rm Ind}_{B(\mathbb F_q)}^{{\rm SL}_2(\mathbb F_q)}
\varepsilon,
\]
where $\varepsilon$ is the quadratic character of $\mathbb F_q^\times$,
extended to $B(\mathbb F_q)$ by
\[
\varepsilon\!\left(
\begin{pmatrix}a&b\\0&a^{-1}\end{pmatrix}
\right)=\varepsilon(a).
\]
To see this directly, $\delta_0$ transforms under $B(\mathbb F_q)$
by $\varepsilon$, and its ${\rm SL}_2(\mathbb F_q)$-span contains $\delta_0$ and 
$
x\longmapsto\psi(bQ_n(x))$ for $b\in\mathbb F_q$. 
These $q+1$ functions are linearly independent and form a basis of
$\mathbb C[V_n]^{G_n}$: the $G_n$-orbits on $V_n$ are $\{0\}$,
the nonzero isotropic vectors, and the $(q-1)$ level sets $Q_n^{-1}(t)$ for $t \neq 0$.
Frobenius reciprocity therefore gives a surjection
\[
{\rm Ind}_{B(\mathbb F_q)}^{{\rm SL}_2(\mathbb F_q)}
\varepsilon
\longrightarrow
\mathbb C[V_n]^{G_n},
\]
which is an isomorphism since both sides have dimension $q+1$.
This induced representation contains no $\sigma\in\mathcal R$ and therefore
$\Theta_{n,\psi}(\sigma)^{G_n}=0$ for $\sigma\in\mathcal R$.

We conclude that
\[
\mathcal A
=
1_{G_n}\oplus
\bigoplus_{\sigma\in\mathcal R}\Theta_{n,\psi}(\sigma)
\]
is a direct summand of both $\mathcal M_\pm$.
By Lemma \ref{lemma-dc-f}, 
\begin{align*}
q
&=|H_n\backslash G_n/K_{n,\alpha}|\\
&=\dim {\rm Hom}_{G_n}(\mathcal M_+,\mathcal M_-)\\
&\geq\dim {\rm End}_{G_n}(\mathcal A)\\
&\geq
1+\sum_{\sigma\in\mathcal R}
\dim {\rm End}_{G_n}\bigl(\Theta_{n,\psi}(\sigma)\bigr)\\
&\geq 1+|\mathcal R|=q.
\end{align*}
Thus we have equality throughout and this forces
\[
\dim {\rm End}_{G_n}\bigl(\Theta_{n,\psi}(\sigma)\bigr)=1
\]
and
\[
{\rm Hom}_{G_n}
\bigl(\Theta_{n,\psi}(\sigma),\Theta_{n,\psi}(\tau)\bigr)=0
\] for $\sigma\neq\tau$.
Thus we have irreducibility and pairwise inequivalence.

Finally, by Frobenius reciprocity, 
\[
q
=
\sum_{\rho \in {\rm Irr}(G_n)}
\dim\rho^{H_n}
\dim\rho^{K_{n,\alpha}}.
\]
The $q$ distinct common constituents $1_{G_n}$ and $\Theta_{n,\psi}(\sigma)$, for $\sigma\in\mathcal R$, 
already contribute at least $q$.
Each must therefore contribute exactly $1$, so both fixed space
dimensions are $1$, and no other $\rho$ can contribute.
\end{proof}

\subsection{Proof of Theorem \ref{thm1}}\label{proofthm1}

For the trivial pair, both sides are $1$, so we consider the nontrivial case. There are exactly $(q-1)$ nontrivial irreducible representations of $G_n$ with fixed vectors for both $H_n$ and $K_{n,\alpha}$ by Lemmas \ref{lemma-dc-f} and \ref{lemma-kriz2}; these are obtained as big theta lifts from ${\rm SL}_2(\mathbb F_q)$ by Lemma \ref{lemma-kriz2}.

Now let $\pi_j=\Theta_{j,\psi}(\sigma)$ and 
$
P_{j,\sigma}\delta_{v^\pm}^j = e_\pm\otimes w_{j,\pm}
$
for $j=1,n$. Let
\[
m_\sigma(s)=\langle e_+,\sigma(s)e_-\rangle \mbox{~and~}
F_j(s,g)=
\left\langle\delta_{v^+}^j,
\omega_{j,\psi}(s,i_j(g))\delta_{v^-}^j\right\rangle,
\]
where $i_1(g)=g$. 
Since $P=P_{j,\sigma}$ is an orthogonal projection commuting with $\omega_{j,\psi}(s,i_j(g))$, we get
\[
F_j(s,g)
=
\left\langle P\delta_{v^+}^j,
\omega_{j,\psi}(s,i_j(g))P\delta_{v^-}^j\right\rangle
+
\left\langle (1-P)\delta_{v^+}^j,
\omega_{j,\psi}(s,i_j(g))(1-P)\delta_{v^-}^j\right\rangle.
\] On the $\sigma$-isotypic component, the action is
\[
\sigma(s)\otimes\pi_j(i_j(g)).
\] The first term is therefore 
\[\left\langle e_+\otimes w_{j,+},
\sigma(s)e_-\otimes\pi_j(i_j(g))w_{j,-}\right\rangle
=
m_\sigma(s)\,
\left\langle w_{j,+},\pi_j(i_j(g))w_{j,-}\right\rangle.
\]
Thus,  
\[
F_j(s,g) = 
m_\sigma(s)
\langle w_{j,+},\pi_j(i_j(g))w_{j,-}\rangle + R_j(s,g)
\]
where for fixed $g$, the function $s\mapsto R_j(s,g)$ is a linear combination of matrix coefficients of irreducible representations other than $\sigma$.
From Schur orthogonality, i.e., 
\[\frac{1}{|G|} \sum_{g \in G} \langle \rho(g)v, v^\prime \rangle \overline{\langle \tau(g)w, w^\prime \rangle} =
\begin{cases}
 \frac{1}{\dim \rho} \langle v, w \rangle \overline{\langle v^\prime, w^\prime \rangle}, &\rho=\tau, \\
 0, &\rho \not\simeq \tau,
\end{cases}
\]
 we get
 \[
\sum_{s\in{\rm SL}_2(\mathbb F_q)}
\overline{m_\sigma(s)}\,R_j(s,g)=0
\]
and 
\[
\frac{1}{|{\rm SL}_2(\mathbb F_q)|}
\sum_{s\in{\rm SL}_2(\mathbb F_q)}
|m_\sigma(s)|^2
=
\frac{\|e_+\|^2\|e_-\|^2}{\dim\sigma}
=\frac{1}{\dim\sigma}.
\]
It follows that
\[
\langle w_{j,+},\pi_j(i_j(g))w_{j,-}\rangle
=
\frac{\dim\sigma}{|{\rm SL}_2(\mathbb F_q)|}
\sum_{s\in {\rm SL}_2(\mathbb F_q)}\overline{m_\sigma(s)} F_j(s,g).
\]
By Lemma \ref{lemma-weil},
\[
F_n(s,g)=q^{-(n-1)}F_1(s,g),
\]
and therefore we get
\[
\langle w_{n,+},\pi_n(i_n(g))w_{n,-}\rangle
=q^{-(n-1)}
\langle w_{1,+},\pi_1(g)w_{1,-}\rangle.
\]
Take $g=g_u^1$, so that $i_j(g)=g_u^j$.
Writing $c_j(u)=c(\pi_j;H_j,K_{j,\alpha,u})$, by Remark \ref{rmk-normsquared}, we get 
\[
c_n(u)
=q^{-2(n-1)}
  \frac{1-q^{-2}}{1-q^{-2n}}\,c_1(u)
=\frac{q^2-1}{q^{2n}-1}\,c_1(u).
\]
The dimension formula then gives
$(\dim\pi_n)c_n(u)=(\dim\pi_1)c_1(u)$.

Corollary \ref{coro-thm1} follows from the rank one result of \cite{ana25}, since $(q^{2n}-1)/(q^2-1)\equiv1\mod p$.

\subsection{A lemma on $\pi_j=\Theta_{j,\psi}(\sigma)$}\label{section-future} 

For future use, we record an elementary lemma about the theta lift. We continue to assume $q$ is odd.

\begin{lemma}\label{lemma-gt-finite}
Let $j \geq 2$. 
Let $\sigma$ be an irreducible complex representation of ${\rm SL}_2(\mathbb F_q)$. Suppose that
\[
\pi_i=\Theta_{i,\psi}(\sigma)
={\rm Hom}_{{\rm SL}_2(\mathbb F_q)}
(\sigma,\omega_{i,\psi})
\]
are nonzero and irreducible for $i=j-1,j$. Let $\Omega_j=M_jN_j$ be the parabolic subgroup of $G_j$ stabilizing $\mathbb F_qe_1$, with
\[
M_j=\{m(a,h)={\rm diag}(a,h,a^{-1})\}
\simeq {\rm GL}_1\times G_{j-1}.
\]
Then there is an $\Omega_j(\mathbb F_q)$-equivariant surjection
\[
\left.\pi_j\right|_{\Omega_j(\mathbb F_q)}
\twoheadrightarrow 1\boxtimes\pi_{j-1},
\]
and thus,
\[
\pi_j\hookrightarrow
{\rm Ind}_{\Omega_j(\mathbb F_q)}^{G_j(\mathbb F_q)}
(1\boxtimes\pi_{j-1}).
\]
\end{lemma}

\begin{proof}
We use the elementary observation that a nonzero invariant functional $\ell:V \to\mathbb C$ gives an equivariant surjection
$
\ell\otimes {\rm id}_W:V\otimes W\twoheadrightarrow W$ by $v\otimes w\longmapsto\ell(v)w$. Write
\[
V_j=\mathbb H\perp V_{j-1},
~
\mathbb H=\mathbb F_qe_1\oplus\mathbb F_qf_1,
~
Q_{\mathbb H}(xe_1+yf_1)=xy.
\]
The tensor product identity, cf. (\ref{tensor}), identifies the Schrödinger model with
\[
\mathbb C[V_j]
=\mathbb C[\mathbb H]\otimes \mathbb C[V_{j-1}],
~
\omega_{j,\psi}
=\omega_{\mathbb H,\psi}\otimes\omega_{j-1,\psi}
\]
as representations of ${\rm SL}_2(\mathbb F_q)$. Define
\[
\ell_{\mathbb H}(\phi)
=\sum_{t\in\mathbb F_q}\phi(te_1).
\]
This functional is ${\rm SL}_2(\mathbb F_q)$-invariant. Indeed, upper unipotents act by multiplication by $\psi(bQ_{\mathbb H})$, which is $1$ on $\mathbb F_qe_1$; diagonal elements only reindex the sum. Let $w = [0, -1; 1, 0]$ be the Weyl element used in \S \ref{theta-finite}. Since the Weil index of $\mathbb H$ is $1$, we get
\[
\ell_{\mathbb H}(\omega_{\mathbb H,\psi}(w)\phi)
=\frac{1}{q}
\sum_{t,x,y\in\mathbb F_q}
\psi(-ty)\phi(xe_1+yf_1)
=\sum_{x\in\mathbb F_q}\phi(xe_1)
=\ell_{\mathbb H}(\phi).
\]
Since $\ell_{\mathbb H}(\delta_0)=1$, we get an ${\rm SL}_2(\mathbb F_q)$-equivariant surjection
\[
\mathcal I_j=\ell_{\mathbb H}\otimes {\rm id}:
\mathbb C[V_j]\twoheadrightarrow\mathbb C[V_{j-1}],
~
(\mathcal I_j\varphi)(v)
=\sum_{t\in\mathbb F_q}\varphi(te_1+v);
\]
indeed it is enough to check this on a pure tensor $\varphi = \phi \otimes \xi$ and then 
for $v\in V_{j-1}$, we have (cf. (\ref{eq-pure-tensor})),
\[
\bigl((\ell_{\mathbb H}\otimes {\rm id})(\phi\otimes \xi)\bigr)(v)
=\left(\sum_{t\in\mathbb F_q}\phi(te_1)\right)\xi(v)\\
=\sum_{t\in\mathbb F_q}\varphi(te_1+v).
\]
We check its equivariance under $\Omega_j(\mathbb F_q)$. For $m(a,h)\in M_j(\mathbb F_q)$,
\[
(\mathcal I_j\omega_{j,\psi}(m(a,h))\varphi)(v)
=\sum_{t\in\mathbb F_q}
\varphi(a^{-1}te_1+h^{-1}v)\\
=(\omega_{j-1,\psi}(h)\mathcal I_j\varphi)(v).
\]
Thus the ${\rm GL}_1(\mathbb F_q)$-factor acts trivially. 
Now for $n\in N_j(\mathbb F_q)$ and $v\in V_{j-1}$, we have
$
n^{-1}e_1=e_1$ and $n^{-1}v-v\in\mathbb F_qe_1$.
Thus, for some $b_n(v)\in\mathbb F_q$,
\[
n^{-1}(te_1+v)=(t+b_n(v))e_1+v.
\]
Reindexing the sum therefore gives
\[
(\mathcal I_j\omega_{j,\psi}(n)\varphi)(v)
=\sum_{t\in\mathbb F_q}
\varphi\bigl((t+b_n(v))e_1+v\bigr)\\
=\sum_{t\in\mathbb F_q}\varphi(te_1+v)
=(\mathcal I_j\varphi)(v).
\]
We have therefore constructed a surjection
\[
\omega_{j,\psi}\twoheadrightarrow
1\boxtimes\omega_{j-1,\psi}
\]
equivariant under
${\rm SL}_2(\mathbb F_q)\times\Omega_j(\mathbb F_q)$. Applying the exact functor
${\rm Hom}_{{\rm SL}_2(\mathbb F_q)}(\sigma,-)$, we get 
\[
\pi_j|_{\Omega_j(\mathbb F_q)}
\twoheadrightarrow 1\boxtimes\pi_{j-1}.
\] By Frobenius reciprocity, we get a nonzero $G_j(\mathbb F_q)$-equivariant map
\[
\pi_j\longrightarrow
{\rm Ind}_{\Omega_j(\mathbb F_q)}^{G_j(\mathbb F_q)}
( 1\boxtimes\pi_{j-1})
\]
which is injective because $\pi_j$ is irreducible by hypothesis.
\end{proof}

\section{Further preliminaries towards a proof of Theorem \ref{thm2}}\label{proof-thm2}

Throughout this section, $F$ is a $p$-adic field of odd residue cardinality $q$. 

\subsection{The pair $(G_n,K_{n,\alpha})$}\label{prelim-thm2}

We collect preliminary results on distinction for the pair $(G_n,K_{n,\alpha})$. This is a multiplicity one pair by \cite{wal12}. Theta correspondence is the key to understand distinction for the pair $(SO(r),SO(r-1))$ and in particular for the pair of our interest which is $(G_n,K_{n,\alpha})$. This is done in \cite[Corollary 11.16]{ggp20} (or see \cite[Proposition 6.1]{gw21}).
\begin{proposition}\label{prop-ggp}
If $\pi$ is an irreducible smooth representation of $G_n$ such that 
\[
{\rm Hom}_{K_{n,\alpha}}(\pi,\mathbb C)\neq 0,
\]
then there exists a unique irreducible genuine smooth representation $\sigma$ of ${\rm Mp}_2(F)$ such that
$\pi\simeq\theta_{n,\psi}(\sigma)$.
\end{proposition}
\begin{proof}
By \cite[Corollary 11.16 (i)]{ggp20}, we have $0 \neq \theta_{n,\psi}^{\leftarrow}(\pi)=\sigma \in {\rm Irr}_{\rm gen}({\rm Mp}_2(F))$. 
The bijection (\ref{eq-howe}) shows that $\sigma$ is the unique irreducible genuine smooth representation of ${\rm Mp}_2(F)$ satisfying $\theta_{n,\psi}(\sigma)\simeq\pi$.
\end{proof}

\begin{remark}\label{atobe}
By \cite[Corollary 11.16 (ii)]{ggp20}, if $\sigma$ is $\psi_\alpha$-generic, then $\pi = \theta_{n,\psi}(\sigma)$ is distinguished by $K_{n,\alpha}$. Indeed, by \cite[Corollary 11.16]{ggp20}, $\pi$ is distinguished by $K_{n,\alpha}$ if and only if $\Theta_{n,\psi}^{\leftarrow}(\pi)$ is $\psi_\alpha$-generic. In particular when $\Theta_{n,\psi}^{\leftarrow}(\pi)$ is irreducible $\psi_\alpha$-genericity follows for $\sigma$ in Proposition \ref{prop-ggp}.  
However, in Proposition \ref{prop-ggp}, $\sigma$ need not be $\psi_\alpha$-generic. Indeed, the small backward theta lift of $1_{SO(3)}$ is the even Weil representation associated with $\psi$ [GS12, Remark following Corollary 9.3, p. 1681], which is not $\psi_\alpha$-generic since $\alpha\notin F^{\times2}$.
\end{remark}

\begin{remark}\label{rmk-tempered}
If $\tau$ is a tempered representation of $G_1$, by Proposition \ref{prop-waldspurger}, we know that 
$\Theta^{\leftarrow}_{1,\psi}(\tau)
\simeq\theta^{\leftarrow}_{1,\psi}(\tau)$, call this $\sigma$. By Remark \ref{atobe}, $\tau$ is $K_{1,\alpha}$-distinguished if and only if $\sigma$ is $\psi_\alpha$-generic. Since $\theta_{1,\psi}(\sigma)\simeq\tau$, persistence (see \S \ref{theta-finite})  implies that $\pi_n=\theta_{n,\psi}(\sigma)\neq 0$ for every $n\geq1$. If $\tau$ is $K_{1,\alpha}$-distinguished, it follows, again by Remark \ref{atobe}, that each $\pi_n$ is $K_{n,\alpha}$-distinguished.
\end{remark}

\subsection{The triple $(G_n,K_{n,\alpha},P_{n,m})$}\label{prelim-thm2-cont}

Next we like to get a better handle on the irreducible genuine smooth representation $\sigma$ of ${\rm Mp}_2(F)$ occurring in Proposition \ref{prop-ggp}. For this, as in the statement of Theorem \ref{thm2}, we are going to assume that the irreducible smooth representation $\pi$ of $G_n$ has a fixed vector for the paramodular subgroup $P_{n,m}$ (cf. \ref{ortho-para}).
We work throughout with the conjugated paramodular subgroups $\mathcal P_{j,c}$ defined in Remark \ref{rmk-para-conj}. Since conjugation preserves the dimensions of fixed-vector spaces, the dimension statements in this and the next section apply equally to these subgroups.

The key input here is the Waldspurger correspondence (cf. Proposition \ref{prop-waldspurger}). The proof of the following result makes crucial use of \cite[Proposition 3.2]{gt16} as well as \cite[Lemma 6.1]{che26}. The same ideas we use in the next section to prove existence and uniqueness of paramodular newforms for Saito Kurokawa representations. We start by recalling the following form of \cite[Lemma 6.1]{che26}, which we need to make use of several times in what follows.

\begin{lemma}[Cheng]\label{lemma-cheng}
Let $\delta$ be an irreducible admissible generic representation of ${\rm GL}_t(F)$, of conductor exponent $a=a(\delta)$, and let $\eta$ be a finite length admissible representation of $G_k$. Write
$
I=\delta\rtimes\eta
$
for normalized parabolic induction to $G_{k+t}$. Suppose that $d\geq0$ satisfies
$
\eta^{\mathcal P_{k,m}}=0$ for $0\leq m<d$. 
Then
$
I^{\mathcal P_{k+t,m}}=0$ for $0\leq m<d+2a$, 
and
\[
\dim I^{\mathcal P_{k+t,d+2a}}
=
\dim\eta^{\mathcal P_{k,d}}.
\]
Moreover, if $\eta$ has no paramodular fixed vectors at any level, neither does $I$.
\end{lemma}

\begin{remark}
Note that there is no genericity assumption on $\eta$. This lemma is stated only for irreducible $\eta$ in \cite[Lemma 6.1]{che26} and there it is assumed that its fixed space at level $d$ has dimension at most one. Neither restriction is needed in Cheng’s fixed vector calculation. Indeed, the calculation in \cite{che26} gives
\[
I^{\mathcal P_{k+t,m}}
\simeq
\bigoplus_{\ell=0}^{\lfloor m/2\rfloor}
\delta^{\Gamma_{t,\ell}}
\otimes
\eta^{\mathcal P_{k,m-2\ell}},
\]
where
\[
\Gamma_{t,\ell}
=
\left\{
g\in{\rm GL}_t(\mathfrak o):
(g_{t1},\ldots,g_{tt})
\equiv(0,\ldots,0,1)\mod{\mathfrak p^\ell}
\right\}.
\]
This decomposition remains valid for finite length $\eta$. Since $\delta^{\Gamma_{t,\ell}}=0$ for $\ell<a$ and $\dim\delta^{\Gamma_{t,a}}=1$, our formulation follows.
\end{remark}

\begin{proposition}\label{prop-split-source}
Let $\sigma$ be an irreducible admissible genuine representation of
${\rm Mp}_2(F)$. Suppose that
\[
\pi=\theta_{n,\psi}(\sigma)\neq 0,~
\pi^{\mathcal P_{n,m}}\neq 0
\]
for some $m\geq 0$. Then the Waldspurger correspondent of $\sigma$
is an irreducible smooth representation of  ${\rm PGL}_2(F)$. Equivalently,
$\theta_{1,\psi}(\sigma)\neq 0$.
\end{proposition}

\begin{proof}
The assertion is immediate when $n=1$, so assume $n\ge2$. By Proposition \ref{prop-waldspurger}, the Waldspurger correspondent of $\sigma$ is an irreducible smooth representation of either ${\rm PGL}_2(F)$ or $PD^\times$. Supposing it is a
representation $\tau_D$ of $PD^\times$ we derive a contradiction.

The representation
$\rho=\theta_{2,\psi}(\sigma)$
is the quaternionic Saito--Kurokawa lift of $\tau_D$ \cite[Proposition 5.8]{sch05}. It is nonzero and
tempered (cf. proof of \cite[Proposition 5.8]{sch05}) and also nongeneric (cf. \cite[Remark on p. 234]{sch05}). 
By \cite[Theorem~7.5.8]{rs07}, a tempered representation of
${\rm PGSp}_4(F)\simeq G_2$ admits paramodular fixed vectors at some level 
if and only if it is generic. It follows that 
$\rho^{\mathcal P_{2,\ell}}=0$ for every $\ell\geq 0$. Thus, $n \neq 2$.

Suppose now that $n \geq 3$. Since $\rho$ is tempered, by \cite[Proposition~3.2]{gt16}, there is
a surjection
\begin{equation}\label{eq-gan-takeda}
I_n =
|\cdot|^{n-\frac{3}{2}}\times
|\cdot|^{n-\frac{5}{2}}\times\cdots\times
|\cdot|^{\frac{3}{2}}\rtimes\rho
\twoheadrightarrow \pi,
\end{equation}
where we take normalized induction. Indeed, in the notation of
that proposition, we take
\[
\dim W=2, \epsilon_0=1,
d=5, d^\prime=2n+1, r=(d^\prime-d)/2 = n-2, \chi_W=1.
\]

The intermediate representations in the successive parabolic inductions defining $I_n$ are all admissible and of finite length. Since $\rho^{\mathcal P_{2,\ell}}=0$ for every $\ell\geq0$, repeated application of the final assertion of Lemma \ref{lemma-cheng}, with $t=1$, gives
$
I_n^{\mathcal P_{n,m}}=0$ for every $m\geq0$.
By the surjection (\ref{eq-gan-takeda}) and exactness of the functor of invariants under a compact open subgroup, we obtain
$
\pi^{\mathcal P_{n,m}}=0$ for every $m\geq0$,
contradicting our hypothesis. Thus it follows that the Waldspurger correspondent of $\sigma$ is an irreducible smooth representation of
${\rm PGL}_2(F)$.
\end{proof}

\begin{remark}\label{rmk-not-yet}
For the first part of Theorem \ref{thm2}, it remains to show that when $\tau=\theta_{1,\psi}(\sigma)$ is tempered, $K_{n,\alpha}$-distinction of $\pi_n$ implies $K_{1,\alpha}$-distinction of $\tau$. This is done in Lemma \ref{lemma-yet}. The nontempered case is treated in the proof of Theorem \ref{thm2}.
\end{remark}

\subsection{Paramodular newforms}\label{section-cheng}

This subsection may be of independent interest. Closely following \cite{che26}, we prove a multiplicity one result for paramodular newforms for representations of Saito Kurokawa type. Indeed \cite[Theorem A (1)]{che26} 
gives analogous vanishing and dimension bounds for generic representations of $SO(n+1,n)$. The proof of Proposition \ref{prop-cheng} follows the fixed vector calculation in the proof of \cite[Lemma 6.1]{che26}, applied successively to induction by unramified characters of ${\rm GL}_1$, starting from Casselman’s newform theorem for $SO(2,1)$ \cite{cas73}. 

\begin{proposition}\label{prop-cheng}
Let $\sigma$ be an irreducible admissible genuine representation of ${\rm Mp}_2(F)$ such that
$
\tau=\theta_{1,\psi}(\sigma)
$
is nonzero and tempered. Then, for $n\ge1$, $\pi_n=\theta_{n,\psi}(\sigma)\neq0$. With $c=a(\tau)$, we have 
$
\pi_n^{P_{n,m}}=0$ for $0\le m<c$, and $\dim\pi_n^{P_{n,c}}\le1$. 
Therefore, if $\pi_n^{P_{n,c}}\neq0$, then the minimal paramodular level of $\pi_n$ is $c$, and its space of paramodular newvectors is one dimensional.
\end{proposition}

\begin{proof}
For $n=1$, this is the well known theorem due to Casselman, so assume $n\ge2$. That $\pi_n \neq 0$ is the persistence principle recalled in \S \ref{theta-finite}. Since $\tau$ is tempered, as in (\ref{eq-gan-takeda}), we have a surjection
\[
\mathcal I_n
=
|\cdot|^{n-\frac{3}{2}}\times
|\cdot|^{n-\frac{5}{2}}\times\cdots\times
|\cdot|^{\frac{1}{2}}\rtimes\tau
\twoheadrightarrow\pi_n.
\]
Set $\mathcal I_1=\tau$ and, recursively,
\[
\mathcal I_j
=
|\cdot|^{j-\frac{3}{2}}\rtimes\mathcal I_{j-1},
\]
for $2\le j\le n$. Note that each $\mathcal I_j$ is admissible and has finite length. By Casselman’s theorem \cite{cas73}, we have $
\tau^{P_{1,m}}=0$ for $0\le m<c$, and $\dim\tau^{P_{1,c}}=1$. 
Therefore, applying Lemma \ref{lemma-cheng} successively with $t=1$ and $d=c$, we obtain
\[
\mathcal I_j^{P_{j,m}}=0 \mbox{~for~} 0\leq m<c \] 
and $\dim\mathcal I_j^{P_{j,c}}=1$, for all $j$. 
From the exactness of invariants under the compact open subgroup $P_{n,m}$, we get a surjection
\[
\mathcal I_n^{P_{n,m}}\twoheadrightarrow
\pi_n^{P_{n,m}}.\] If $\pi_n^{P_{n,c}}\neq0$, the bound is an equality, and we get minimality from the vanishing at all smaller levels.
\end{proof}

Retain the hypotheses and notation of Proposition \ref{prop-cheng}, and put $\pi_j=\theta_{j,\psi}(\sigma)$ for $1\le j\le n$. We now construct the paramodular vectors recursively, starting from the rank one newvector. Suppose that $c=a(\tau)=2r>0$.

For $j\geq 2$, let $\Omega_j$ be the parabolic subgroup of $G_j$ stabilizing $Fe_1$. Write $\Omega_j=M_jN_j$, where
$
m_j(a,h)={\rm diag}(a,h,a^{-1})\in M_j
\simeq {\rm GL}_1(F)\times G_{j-1}.
$
Since $\mathcal I_j=|\cdot|^{j-\frac{3}{2}}\rtimes\mathcal I_{j-1}$ is normalized induction, its sections satisfy
\[
F(m_j(a,h)ng)
=
|a|^{2j-2}\mathcal I_{j-1}(h)F(g).
\]
Indeed, $\delta_{\Omega_j}(m_j(a,h))=|a|^{2j-1}$, and hence
$
\delta_{\Omega_j}^{1/2}(m_j(a,h))\,|a|^{j-\frac{3}{2}}
=
|a|^{2j-2}.
$
By \cite[Corollaries 3.3 and 5.2]{che26}, we have
\[
{\rm pr}_{M_j}(\Omega_j\cap\mathcal P_{j,c})
=
\mathfrak o^\times\times\mathcal P_{j-1,c}.
\]
\begin{remark}
Here, Cheng’s $K^0_{j,c}$ is our $\mathcal P_{j,c}$. We apply \cite[Corollary 5.2]{che26} to the ${\rm GL}_1$-Levi, with $m=c$ and $i=c/2$ to get
${\rm pr}_{M_j} \bigl(\Omega_j\cap s\mathcal P_{j,c} s^{-1}\bigr) = \mathfrak o^\times\times\mathcal P_{j-1,c}$ for an element $s$ which we know belongs to $\mathcal P_{j,c}$ by \cite[Corollary 3.3]{che26}.
\end{remark}
Now choose a nonzero rank one newvector
$
v_{1,c}\in\tau^{\mathcal P_{1,c}},
$
and set $F_{1,c}=v_{1,c}$. For $2\leq j\leq n$, define $F_{j,c}\in\mathcal I_j$ recursively to be right $\mathcal P_{j,c}$-invariant, supported on $\Omega_j\mathcal P_{j,c}$, and normalized by
$
F_{j,c}(1)=F_{j-1,c}.
$
Explicitly,
\begin{equation}\label{eq-paranormal}
F_{j,c}(m_j(a,h)nk)
=
|a|^{2j-2}\mathcal I_{j-1}(h)F_{j-1,c},
\end{equation}
for $k\in\mathcal P_{j,c}$ 
and $F_{j,c}$ is zero outside $\Omega_j\mathcal P_{j,c}$.

Note that this definition is independent of the chosen decomposition as the ${\rm GL}_1$-component of an element of $\Omega_j\cap\mathcal P_{j,c}$ is a unit, and also $F_{j-1,c}$ is fixed by its $G_{j-1}$-component. Since $F_{j,c}(1)\neq0$, the section is nonzero, and from the proof of Proposition \ref{prop-cheng} we see that $\mathcal I_j^{\mathcal P_{j,c}}=\mathbb C F_{j,c}$.

By \cite[Proposition 3.2]{gt16}, $\pi_j$ is the unique irreducible quotient of the standard module $\mathcal I_j$. By the exactness of parabolic induction we get a surjection
\[
\mathcal I_j\twoheadrightarrow
I_j^+=|\cdot|^{j-\frac{3}{2}}\rtimes\pi_{j-1}.
\]
Since $I_j^+$ is nonzero and of finite length, it has an irreducible quotient, which then must be $\pi_j$.

Now we fix the quotient maps
\[
p_j:I_j^+
=
|\cdot|^{j-\frac{3}{2}}\rtimes\pi_{j-1}
\twoheadrightarrow\pi_j
\]
and define 
\[
q_j:\mathcal I_j\twoheadrightarrow\pi_j,\]
\begin{equation}\label{eq-compatible}
q_1={\rm id}_\tau,
~
q_j(F)=p_j\bigl(g\longmapsto q_{j-1}(F(g))\bigr).
\end{equation}
Now the newform in $\pi_j$ is defined as
\[
v_{j,c}=q_j(F_{j,c})\in\pi_j^{\mathcal P_{j,c}}.
\]
Thus, in the induced model $I_j^+$, the section defining $v_{j,c}$ takes the value $v_{j-1,c}$ at the identity. 

\begin{remark}\label{vjc-nonvanishing}
Although $F_{j,c}\neq0$, its image $v_{j,c}$ in the irreducible quotient may a priori vanish. We still need to prove its nonvanishing. Once $v_{j,c}\neq0$ is proved, the preceding proposition shows that
\[
\pi_j^{\mathcal P_{j,c}}=\mathbb C v_{j,c}
\]
and that $c$ is the minimal paramodular level of $\pi_j$. For better clarity, we postpone this to Proposition \ref{prop-rsj} in the next section.
\end{remark}

\subsection{Nonvanishing of the newvectors}

We prove that the vectors $F_{j,c}$ have nonzero images in the irreducible quotients $\pi_j$. Towards this, first we collect a couple of results from \cite{jan96,rs07} as a lemma.

For a representation $\rho$ of $G_1$, let
\[
\mathcal I_j^-(\rho)
=
|\cdot|^{-j+\frac{3}{2}}\times\cdots\times|\cdot|^{-\frac{1}{2}}\rtimes\rho,\]
with $\mathcal I_1^-(\rho)=\rho$. 
For tempered $\tau$, note that the representation $\pi_j ={\rm LQ}(\mathcal I_j)$ is the unique irreducible subrepresentation of $\mathcal I_j^-(\tau)$.

\begin{lemma}[Jantzen, Roberts–Schmidt]\label{lemma-jantzen}
Let $\tau$ be an irreducible tempered representation of $G_1$, and let $c=a(\tau)$.
\begin{enumerate}
\item\label{i}
If $\tau$ is supercuspidal, then, for $j\geq2$, there is an exact sequence
\[
0\rightarrow\pi_j
\rightarrow|\cdot|^{-j+\frac{3}{2}}\rtimes\pi_{j-1}
\rightarrow\eta_j
\rightarrow 0,
\]
where $\eta_j$ is irreducible. Moreover,
$
\eta_2^{\mathcal P_{2,c}}=0,
$
and, for $j\geq3$,
\[
\eta_j\hookrightarrow
|\cdot|^{2-j}{\rm St}_{{\rm GL}_2}\rtimes
\mathcal I_{j-2}^-(\tau).
\]
\item\label{ii} 
For $j \geq 2$, let $\mu$ be a ramified character of $F^\times$. The irreducible constituents of
$
\mu\rtimes 1_{G_{j-1}}
$
occur with multiplicity one and correspond bijectively to the irreducible constituents $\rho$ of  ${\rm Ps}(\mu,\mu^{-1})$, viewed as a representation of $G_1$. The constituent corresponding to $\rho$ embeds into $\mathcal I_j^-(\rho)$; when $\rho=\tau$, it is $\pi_j$.
\end{enumerate}
\end{lemma}

\begin{proof}
Under $G_2\simeq{\rm PGSp}_4$, our $|\cdot|^s\rtimes\tau$ corresponds to $|\cdot|^s\tau\rtimes|\cdot|^{-s}$ in \cite{rs07}. Thus, for (\ref{i}), in the case $j=2$, since $\tau$ is supercuspidal, by \cite[Group XI on p. 40 and Proposition 5.5.1 (iv)]{rs07}, we know that $\pi_2$ is irreducible nontempered and nongeneric whereas $\eta_2$ is irreducible essentially square integrable and generic. By \cite[Proposition 5.5.5]{rs07}, we know that 
$\pi_2$ has minimal paramodular level $c$ and $\eta_2$ has minimal paramodular level $c+1$. 
For $j\geq3$, we apply \cite[Proposition 3.1 (2)]{jan96} with
\[
\rho=1, \sigma=\tau,
\ell=j-2,\alpha=-\ell-\tfrac12.
\]
Its reducibility hypothesis, namely \cite[C1/2 on p.6]{jan96}, follows from \cite[Groups X and XI on p. 40]{rs07}. In the notation of \cite[p. 1]{jan96}, we have $\zeta(1,j-2;\tau)=\pi_{j-1}$ and therefore $\pi = |\cdot|^{-j+\frac{3}{2}}\rtimes\pi_{j-1}$ in \cite[Proposition 3.1]{jan96}. By \cite[Proposition 3.1 (2)]{jan96}, the representation $\pi$ has two distinct irreducible subquotients given by 
\[
\Pi_j = L([|\cdot|^{-j+\frac{3}{2}},|\cdot|^{-\frac{1}{2}}];\tau),
\]
and
\[
\eta_j
=
L\!\left(
|\cdot|^{2-j}{\rm St}_{{\rm GL}_2},
[|\cdot|^{-j+\frac{7}{2}},|\cdot|^{-\frac{1}{2}}];\tau
\right),
\]
in the notation of \cite[p. 5]{jan96}.
Thus, we get
\[
\eta_j\hookrightarrow
|\cdot|^{2-j}{\rm St}_{{\rm GL}_2}
\times|\cdot|^{-j+\frac{7}{2}}\times\cdots\times|\cdot|^{-\frac{1}{2}}
\rtimes\tau
=
|\cdot|^{2-j}{\rm St}_{{\rm GL}_2}\rtimes
\mathcal I_{j-2}^-(\tau).
\]
By the exactness of parabolic induction,
\[
\pi= |\cdot|^{-j+\frac{3}{2}}\rtimes\pi_{j-1}
\hookrightarrow
|\cdot|^{-j+\frac{3}{2}}\rtimes\mathcal I_{j-1}^-(\tau)
=\mathcal I_j^-(\tau).
\]
Since $\pi$ is nonzero and of finite length, it contains an irreducible subrepresentation, necessarily $\Pi_j=\pi_j$, and thus it follows that we have
\[
0\rightarrow\pi_j
\rightarrow|\cdot|^{-j+\frac{3}{2}}\rtimes\pi_{j-1}
\rightarrow\eta_j
\rightarrow 0.
\]
For (\ref{ii}), we need to follow \cite[\S 5]{jan96}, where the representation 
$\pi = |\cdot|^\alpha \zeta(\rho_0, k) \rtimes \zeta (\rho,\ell;\sigma)$ for $\rho_0 \not\simeq \rho$ is studied. We write $\mu=|\cdot|^s\mu_0$ with $\mu_0$ unitary, and apply \cite[Lemma 5.4 (1)]{jan96} with
\[
\rho_0=\mu_0, \rho=1,
\sigma=1_{G_0}, k=1,
\ell=j-1, \alpha=s.
\]
Ramification ensures $\mu_0\neq1$. We have $\pi = \mu\rtimes1_{G_{j-1}}$. Also (cf. \cite[Proposition 5.3]{jan96}),
\[
|\cdot|^\alpha\zeta(\rho_0,k)\rtimes\sigma
=\mu\rtimes1_{G_0}
\simeq\mu\times\mu^{-1},
\]
under $G_1\simeq{\rm PGL}_2(F)$. With $L(\Delta_i)$ being its irreducible constituents, and thus $i$ is either $1$ or $2$, we get  
\[
\Pi_i=L\bigl([|\cdot|^{-j+\frac{3}{2}},|\cdot|^{-\frac{1}{2}}],\Delta_i\bigr)
\]
as the irreducible constituents of $\pi$ by \cite[Proposition 5.3]{jan96}, and by \cite[Lemma 5.4 (1)]{jan96}, we have 
\[
\Pi_i
\hookrightarrow
|\cdot|^{-\frac{j-1}{2}}\zeta(1,j-1)\rtimes L(\Delta_i)
\hookrightarrow
|\cdot|^{-j+\frac{3}{2}}\times\cdots\times|\cdot|^{-\frac{1}{2}}\rtimes L(\Delta_i)
=\mathcal I_j^-(L(\Delta_i)).
\]
When $L(\Delta_i)=\tau$, the uniqueness of the irreducible subrepresentation of $\mathcal I_j^-(\tau)$ identifies $\Pi_i$ with $\pi_j$.
\end{proof}

The following proposition is the result that we are after.
\begin{proposition}\label{prop-rsj}
Assume that $\tau$ is tempered and $c=a(\tau)=2r>0$. Then the recursively defined vectors
\[
v_{j,c}=q_j(F_{j,c})
\]
are nonzero for $1\leq j\leq n$. Hence,
\[
\pi_j^{\mathcal P_{j,c}}=\mathbb C v_{j,c},
\]
and $v_{j,c}$ is a paramodular newform of minimal level $c$. Moreover, $F_{j,c}$ generates $\mathcal I_j$ as a $G_j$-representation.
\end{proposition}

\begin{proof}
Recall that, by Proposition \ref{prop-cheng},
$\pi_j^{\mathcal P_{j,m}}=0$ for $m < c$ and $
\mathcal I_j^{\mathcal P_{j,c}}=\mathbb C F_{j,c}$. We apply Lemma \ref{lemma-cheng} to get 
\[
\mathcal I_s^-(\tau)^{\mathcal P_{s,m}}=0\]
for $m<c$.
We prove that $\dim\pi_j^{\mathcal P_{j,c}}=1$. The case $j=1$ is the classical newform theorem of Casselman \cite{cas73}.

Suppose first that $\tau$ is supercuspidal. By Lemma \ref{lemma-jantzen} (\ref{i}), $\eta_2^{\mathcal P_{2,c}}=0$. For $j\geq3$, now apply Lemma \ref{lemma-cheng} with $t=2$ and we conclude
\[
\eta_j^{\mathcal P_{j,c}}=0,
\]
since $a(|\cdot|^{2-j}{\rm St}_{{\rm GL}_2})=1$. Taking invariants in the exact sequence and applying Lemma \ref{lemma-cheng} with $t=1$, we obtain inductively
\[
\dim\pi_j^{\mathcal P_{j,c}}
=
\dim\pi_{j-1}^{\mathcal P_{j-1,c}}
=1.
\]
Suppose now that $\tau$ is not supercuspidal. Since $c$ is positive and even, we may choose a ramified character $\mu$ with $2a(\mu)=c$ such that $\tau$ is a constituent of ${\rm Ps}(\mu,\mu^{-1})$. Explicitly, take $\mu=\chi$ if $\tau={\rm Ps}(\chi,\chi^{-1})$, and $\mu=\xi|\cdot|^{-\frac{1}{2}}$ if $\tau={\rm St}_{{\rm GL}_2} \otimes \xi$, where $\xi$ is ramified quadratic. 
In the Steinberg case, the only other constituent is the ramified character  $\xi\circ\det$, which has no paramodular fixed vectors at any level. Therefore, $\mathcal I_j^-(\xi\circ\det)$ does not admit a paramodular fixed vector by Lemma \ref{lemma-cheng}. Thus, by Lemma \ref{lemma-jantzen} (\ref{ii}), any constituent $\rho$ of $\mu\rtimes 1_{G_{j-1}}$ other than $\pi_j$ will have no paramodular fixed vectors at any level.

Now if we apply Lemma \ref{lemma-cheng} with $\delta=\mu$, $\eta=1_{G_{j-1}}$, and $d=0$, we get 
\[
\dim
\bigl(\mu\rtimes 1_{G_{j-1}}\bigr)^{\mathcal P_{j,c}}
=1.
\]
From the exactness of compact open invariants it follows that 
\[
\dim\pi_j^{\mathcal P_{j,c}}=1.
\]
Thus we get an isomorphism
\[
q_j:\mathcal I_j^{\mathcal P_{j,c}}
\xrightarrow{\ \sim\ }
\pi_j^{\mathcal P_{j,c}},
\]
and therefore $v_{j,c}=q_j(F_{j,c})\neq0$.

Since $\mathcal I_j$ is a standard module, we know that $\ker q_j$ is its unique maximal proper subrepresentation. Every proper subrepresentation is contained in $\ker q_j$, whereas $F_{j,c}\notin\ker q_j$ and therefore it follows that $F_{j,c}$ generates $\mathcal I_j$.
\end{proof}

\subsection{Proof of Theorem \ref{thm-cheng}}\label{section-crsj}

If $c > 0$, Theorem \ref{thm-cheng} is immediate from Proposition \ref{prop-cheng} and Proposition \ref{prop-rsj}. If $c=0$, then $\tau$ is unramified, and the Langlands-quotient description in \cite[Proposition 3.2]{gt16} shows that $\pi_n$ is unramified. 
Since $P_{n,0}=G_n(\mathfrak o)$, we have $\pi_n^{P_{n,0}}\ne0$, and Proposition \ref{prop-cheng} gives $\dim\pi_n^{P_{n,0}}=1$.

\subsection{Mackey theory}\label{mackey}

We follow the argument in \cite[Proof of Theorem 15.1]{ggp12} simplified to our situation. Recall that, for $j\geq 2$, $\Omega_j$ is the parabolic subgroup of $G_j$ stabilizing $Fe_1$, with Levi subgroup ${\rm GL}_1(F)\times G_{j-1}$. For a nontrivial irreducible admissible representation $\rho$ of $G_{j-1}$, let
\[
I_j(\rho)
=
{\rm Ind}_{\Omega_j}^{G_j}
\bigl(|\cdot |^{j-\frac{3}{2}}\boxtimes \rho \bigr).
\]
Recall that $I_j(\rho)$ is realized as the space of smooth functions $F: G_j \to \rho$ satisfying
\begin{equation}\label{eq-induction}
F(m(a,h)ng)
=
\delta_{\Omega_j}^{1/2}(m(a,h))|a|^{j-\frac{3}{2}}
\rho(h)F(g)
=
|a|^{2j-2}\rho(h)F(g),
\end{equation}
where $m(a,h)={\rm diag}(a,h,a^{-1})$, $n$ belongs to the unipotent radical of $\Omega_j$, and $g\in G_j$,  with $G_j$ acting by right translation.

Note that the action of $K_{j,\alpha}$ on $\Omega_j\backslash G_j$ has two orbits. The closed orbit consists of the isotropic lines contained in $(e_j+\alpha f_j)^\perp$ and has stabilizer $R_{j,\alpha}=\Omega_j\cap K_{j,\alpha}$. The open orbit consists of the isotropic lines not contained in $(e_j+\alpha f_j)^\perp$ and is represented by $Fe_j$ and has stabilizer
$
S_j= SO \bigl((Fe_j\oplus Ff_j)^\perp\bigr)
\simeq G_{j-1}.
$
Let $w_j\in G_j$ interchange $e_1$ with $e_j$ and $f_1$ with $f_j$, fixing the remaining basis vectors. Then $\Omega_jw_j$ represents the open orbit, and conjugation by $w_j$ identifies $S_j$ with the $G_{j-1}$-factor of the Levi subgroup of $\Omega_j$. 
By Mackey theory about restriction of an induced representation to a
subgroup, we get
\[
0\longrightarrow
{\rm ind}_{S_j}^{K_{j,\alpha}}\rho^{w_j}
\rightarrow
I_j(\rho)\big|_{K_{j,\alpha}}
\xrightarrow{ {\rm res}}
{\rm Ind}_{R_{j,\alpha}}^{K_{j,\alpha}}
\left(
|\cdot |^{j-1}\boxtimes
\rho \big|_{K_{j-1,\alpha}}
\right)
\rightarrow 0,
\]
where ${\rm res}(F)=F|_{K_{j,\alpha}}$ and $\rho^{w_j}(s)=\rho(w_jsw_j^{-1})$.
Applying ${\rm Hom}_{K{j,\alpha}}(-,\mathbb C)$, by Frobenius reciprocity \cite[Corollary 2.4.3]{cas95}, we have 
\[
0\rightarrow
{\rm Hom}_{K_{j-1,\alpha}}(\rho,\mathbb C)
\rightarrow
{\rm Hom}_{K_{j,\alpha}}(I_j(\rho),\mathbb C)
\rightarrow
{\rm Hom}_{G_{j-1}}(\rho,\mathbb C).
\]
The last space vanishes as $\rho$ is irreducible and nontrivial. Thus,
\begin{equation}\label{eq-iso}
{\rm Hom}_{K_{j,\alpha}}(I_j(\rho),\mathbb C)
\simeq
{\rm Hom}_{K_{j-1,\alpha}}(\rho,\mathbb C).
\end{equation}

We can make the identification in (\ref{eq-iso}) precise as follows. For a functional $\lambda$ on the right, define
\[
L_\lambda(F) = 
\int_{C_{j,\alpha}}\lambda(F(c))\,dc,
\]
where $C_{j,\alpha}=K_{j,\alpha}\cap G_j(\mathfrak o)$. Indeed, the above integral converges absolutely by compactness of $C_{j,\alpha}$ and local constancy of the integrand. It defines a $K_{j,\alpha}$-invariant functional on $I_j(\rho)$ because we have, by the Iwasawa decomposition, $K_{j,\alpha}=R_{j,\alpha}C_{j,\alpha}$, and 
the function $f_F(h)=\lambda(F(h))$
satisfies, for $r=m_j(a,k)n\in R_{j,\alpha}$ and $h\in K_{j,\alpha}$,
\[
f_F(rh)
=\lambda\!\left(F(m_j(a,k)nh)\right)
=|a|^{2j-2}\lambda\!\left(\rho(k)F(h)\right)
=|a|^{2j-2}\lambda(F(h))
=\delta_{R_{j,\alpha}}(r)f_F(h),
\]
by (\ref{eq-induction}) and the $K_{j-1,\alpha}$-invariance of $\lambda$ (cf. \cite[Proposition 3.1.3 and its proof]{cas95}). Moreover, if $\lambda$ is nonzero then $L_\lambda$ is also nonzero; to see this, choose $v$ with $\lambda(v)=1$ and a closed orbit section supported on $R_{j,\alpha}J$, taking the value $v$ on a sufficiently small compact open subgroup $J\subset C_{j,\alpha}$; then lift this section to $I_j(\rho)$ using the surjectivity of restriction.

We have thus proved the following lemma.
\begin{lemma}\label{lemma-mackey}
Let $\rho$ be a nontrivial irreducible admissible representation of $G_{j-1}$, and let 
\[
I_j(\rho)=
{\rm Ind}_{\Omega_j}^{G_j}
\bigl(|\cdot|^{j-\frac{3}{2}}\boxtimes\rho\bigr).
\]
Then we have an isomorphism
\[
{\rm Hom}_{K_{j-1,\alpha}}(\rho,\mathbb C)
\xrightarrow{\;\sim\;}
{\rm Hom}_{K_{j,\alpha}}(I_j(\rho),\mathbb C),
\]
given by $\lambda\mapsto L_\lambda$.
\end{lemma}

\subsection{Proof of Remark \ref{rmk-not-yet}}\label{section-yet}

Thanks to the Mackey theory argument in \S \ref{mackey} we can now complete the unfinished task in Remark \ref{rmk-not-yet}. We state this as a lemma.

\begin{lemma}\label{lemma-yet}
Let $\pi_j=\theta_{j,\psi}(\sigma)$, where
$\sigma$ is an irreducible admissible genuine representation of ${\rm Mp}_2(F)$, 
and suppose that $\tau=\pi_1$ is nonzero and tempered. If $\pi_n$ is $K_{n,\alpha}$-distinguished, then $\pi_j$ is $K_{j,\alpha}$-distinguished for every $1\leq j\leq n$. In particular, $\tau$ is $K_{1,\alpha}$-distinguished.
\end{lemma}

\begin{proof}
Since $\pi_1$ is tempered, by \cite[Proposition 3.2]{gt16} and by the exactness and transitivity of parabolic induction, we get a surjection
\[
p_j:I_j^+ 
=
{\rm Ind}_{\Omega_j}^{G_j}
\bigl(|\cdot |^{j-\frac{3}{2}}\boxtimes \pi_{j-1} \bigr)
\twoheadrightarrow\pi_j.
\]
This gives rise to an injection
\begin{equation}\label{eq-injection}
p_j^*:
{\rm Hom}_{K_{j,\alpha}}(\pi_j,\mathbb C)
\hookrightarrow
{\rm Hom}_{K_{j,\alpha}}(I_j^+,\mathbb C) \simeq {\rm Hom}_{K_{j-1,\alpha}}(\pi_{j-1},\mathbb C),
\end{equation}
where the second identification is by Lemma \ref{lemma-mackey}.
Starting from $j=n$, descending induction shows that
${\rm Hom}_{K_{j,\alpha}}(\pi_j,\mathbb C)\ne0$
for every $1\le j\le n$.
\end{proof}

\subsection{$K_{n,\alpha}$-invariant form}\label{invariant}

We continue with the assumptions of \S \ref{section-yet}. Thus for each $j \geq 1$, the representation
$\pi_j=\theta_{j,\psi}(\sigma)$ is nonzero and $K_{j,\alpha}$-distinguished. Moreover, $\pi_1$ is tempered. To define a nontrivial $K_{n,\alpha}$-invariant linear form on $\pi_n$, we start with a nonzero toric functional
\[
\Lambda_{1,\alpha}\in {\rm Hom}_{K_{1,\alpha}}(\pi_1,\mathbb C).
\]
Since $K_{1,\alpha}\simeq E^\times/F^\times$ is compact, this linear form may be obtained by averaging a suitable linear functional over $K_{1,\alpha}$. 
Consider the functionals $L_{\Lambda_{j-1,\alpha}}$ given by
\[
L_{\Lambda_{j-1,\alpha}}(F)
=
\int_{C_{j,\alpha}}
\Lambda_{j-1,\alpha}(F(c))dc,
\]
where we take ${\rm vol}(C_{j,\alpha})=1$. Under our hypotheses, the map $p_j^*$ in (\ref{eq-injection}) is an isomorphism (of one dimensional spaces).
So we get a unique form $\Lambda_{j,\alpha}$ satisfying
\[
\Lambda_{j,\alpha}\circ p_j
=
L_{\Lambda_{j-1,\alpha}},
\]
or explicitly,
\begin{equation}\label{eq-recursive}
\Lambda_{j,\alpha}(p_j(F))
=
\int_{C_{j,\alpha}}
\Lambda_{j-1,\alpha}(F(c))\,dc.
\end{equation}
Note that its normalization is fixed by the initial toric functional, the quotient maps, and the chosen  Haar measures.
For $g\in G_1$, let $h=i_n(g)$ and $K_{n,\alpha,h}=hK_{n,\alpha}h^{-1}$. The translated functional
\begin{equation}\label{eq-translated}
\Lambda_{n,\alpha,h}(w)
=
\Lambda_{n,\alpha}\bigl(\pi_n(h^{-1})w\bigr)
\end{equation}
is then $K_{n,\alpha,h}$-invariant.

\begin{remark}\label{rmk-k-def}
For $u\in\mathfrak o$, write
$\Lambda_{n,\alpha,u}=\Lambda_{n,\alpha,g_u}$, where
$g_u=i_n(g_u^1)$ as in (\ref{k-def}). This functional is $K_{n,\alpha,u}$-invariant.
\end{remark}

\section{Proof of Theorem \ref{thm2}}\label{proofthm2added}

In this section, we prove Theorem \ref{thm2} in \S \ref{sub-pf-thm2}. The idea is to reduce the period evaluation on newform one rank at a time from $j+1$ to $j$. Towards this we make use of the recursive constructions, both of the period and of the newform, of \S \ref{proof-thm2}. We show that these periods differ by the volume of a suitable double coset intersection. This step is reminiscent of the role played by Lemma \ref{lemma-weil} in \S \ref{proof-thm1}. We compute this volume by counting isotropic lines modulo $\mathfrak p^r$ in two quadratic lattices that we introduce in \S \ref{section-ql}. The necessary preliminary computations are carried out in \S \ref{section-lemma-dc}. 

\subsection{Quadratic lattices}\label{section-ql}

Recall our notations from \S \ref{ortho-para}; we have $\mathfrak o=\mathfrak o_F$, $\alpha\in \mathfrak o^\times\setminus \mathfrak o^{\times 2}$, and $u\in\mathfrak o$. We have $V_n=U_{n-1}\perp T_n$, where $T_n=\langle e_n,v_0,f_n\rangle$ (cf. (\ref{directsum})). The subgroup $K_{n,\alpha,u}=g_uK_{n,\alpha}g_u^{-1}$ (cf. (\ref{k-def})) is the stabilizer of 
\[
\xi_{\alpha,u}
=g_u(e_n+\alpha f_n)
=(1-\alpha u^2)e_n+\alpha u\,v_0+\alpha f_n,
\]
and $H_n ={\rm Stab}_{G_n}(v_0)$. With 
\[L_{n,m}= \mathfrak o e_1 \oplus \dots \oplus \mathfrak o e_n \oplus \mathfrak p^m v_0 \oplus \mathfrak p^m f_n \oplus \dots \oplus \mathfrak p^m f_1 \]
as in \S \ref{ortho-para}, we consider the lattices
$
L_X = L_{n,0}\cap\xi_{\alpha,u}^{\perp}$ and $
L_Y = L_{n,0}\cap\langle\xi_{\alpha,u},v_0\rangle^{\perp}.
$
Note that $1-\alpha u^2\in\mathfrak o^\times$. The three vectors
\[
\xi_{\alpha,u},~
v_0-u\xi_{\alpha,u},~
(1-\alpha u^2)e_n-\alpha f_n
\]
form an orthogonal $\mathfrak o$-basis of $L_{n,0}\cap T_n$, with respective norms
$
\alpha$, $1-\alpha u^2$ and $-\alpha(1-\alpha u^2)$, all of which are units. We have,
\begin{align*}
L_{n,0} &=\mathfrak o\xi_{\alpha,u}
\perp L_X, \\
L_X &= L_Y 
\perp\mathfrak o(v_0-u\xi_{\alpha,u}),\\
L_Y &=(L_{n,0}\cap U_{n-1})
\perp
\mathfrak o\bigl((1-\alpha u^2)e_n-\alpha f_n\bigr).
\end{align*}
Both $L_X$ and $L_Y$ are unimodular. Their mod $\mathfrak p$ reductions have dimensions $2n$ and $2n-1$ and Witt index $n-1$. 

The reduction of $L_X$ is nonsplit as its binary part 
\[\langle v_0-u\xi_{\alpha,u}, (1-\alpha u^2)e_n-\alpha f_n \rangle
 \]
has quadratic form $\overline{1-\alpha u^2}x^2+\overline{-\alpha(1-\alpha u^2)}y^2 =
(1-\bar\alpha\bar u^2)(x^2-\bar\alpha y^2)$, 
which is anisotropic over $  F_q$. For $L_Y$, the reduction of the one dimensional part is anisotropic; indeed, the remaining summand is the line
$
\mathfrak o\bigl((1-\alpha u^2)e_n-\alpha f_n\bigr)
$,
whose reduction has quadratic form
$
-\bar\alpha(1-\bar\alpha\bar u^2)y^2
$, and as the coefficient is nonzero, this vanishes only when $y=0$. 

Recall from \S \ref{mackey} that $C_{n,\alpha} = K_{n,\alpha}\cap G_n(\mathfrak o)$.
Since $g_u\in G_n(\mathfrak o)$, we have $C = g_uC_{n,\alpha}g_u^{-1}
=K_{n,\alpha,u}\cap G_n(\mathfrak o)$; this is precisely the integral isometries fixing $\xi_{\alpha,u}$. Therefore,
\[
C \cap H_n(\mathfrak o)
=\{h\in G_n(\mathfrak o) \mid h\xi_{\alpha,u}=\xi_{\alpha,u}, ~ hv_0=v_0\},
\]
and in the notation of Remark \ref{rmk-para-conj},
\[
C \cap C_n(r)
=
\{h\in C_n(r) \mid h\xi_{\alpha,u}=\xi_{\alpha,u}\}.
\]
We have the identifications
\[C = SO(L_X), ~ C \cap H_n(\mathfrak o) = SO(L_Y). \]

\subsection{A lemma on double cosets}\label{section-lemma-dc}

Fix $n \geq 2$ and $r\geq 1$. Let $C=g_uC_{n,\alpha}g_u^{-1}$ and $\ell_0=Fe_1$. Recall that $\Omega_n= {\rm Stab}_{G_n}(\ell_0)$, a maximal parabolic in $G_n$. We are interested in the intersection of $C$ with the identity double coset for $(\Omega_n,\mathcal P_{n,2r})$ and we claim that this intersection is precisely the identity double coset for the corresponding subgroups $(C \cap \Omega_n,C \cap \mathcal P_{n,2r})$. 

Let $X_r$ and $Y_r$ denote the sets of isotropic lines, i.e., rank one direct summands on which the quadratic form vanishes, in the reductions modulo $\mathfrak p^r$ of the two lattices $L_X$ and $L_Y$ that we considered in \S \ref{section-ql}. The integral orthogonal decompositions in \S \ref{section-ql} identify $Y_r$ with a subset of $X_r$. Define
$
\rho_r:C\rightarrow X_r$ by $
c \mapsto\overline{c^{-1}\mathfrak o e_1}$. Then
\[
\rho_r(c)\in Y_r
\iff  
B_n(c^{-1}e_1,v_0)\in\mathfrak p^r.
\]
We have,
\begin{lemma}\label{lemma-dc}
The group $C$ (resp. $C\cap H_n(\mathfrak o)$) acts transitively on $X_r$ (resp. $Y_r$). Every line $\bar{\ell}$ in $X_r$ has a lift to an isotropic line
$
\mathfrak o v\subset L_{n,0}\cap\xi_{\alpha,u}^{\perp}
$
with $v$ being primitive (i.e., $v\notin\mathfrak pL_{n,0}$). Moreover, any two lifts of the same line are carried to one another by an element of $C\cap C_n(r)$. We also have the following double coset identity:
\[
C\cap \Omega_n\mathcal P_{n,2r}
=\rho_r^{-1}(Y_r)
=(C\cap \Omega_n)(C\cap\mathcal P_{n,2r}).
\]
\end{lemma}

\begin{proof}
Let $L$ denote either $L_X$ or $L_Y$ so that $SO(L)$ is respectively $C$ or $C \cap H_n(\mathfrak o)$ (cf. \S \ref{section-ql}).
We will always regard an integral isometry of either lattice as an element of $G_n(\mathfrak o)$ by extending it by the identity on its orthogonal complement in $L_{n,0}$. 

For the proof, we work with the Eichler transvections \cite[\S 3.1]{ghs09}. Explicitly, for $t\in L\cap\langle e_1,f_1\rangle^\perp$, we define
\[
T_t(e_1)=e_1+t-Q(t)f_1,~
T_t(f_1)=f_1,~
T_t(z)=z-B_n(z,t)f_1,\]
for $z\in L\cap\langle e_1,f_1\rangle^\perp$. 
Using the orthogonal decomposition and the polarization identity
\[
Q(z+xt)=Q(z)+xB_n(z,t)+x^2Q(t),
\]
we see that $Q(T_t(v))=Q(v)$.
Its matrix is integral and triangular with all diagonal entries equal to $1$. Hence $T_t\in SO(L)$.
Moreover,
\[
t\in\mathfrak p^rL
\implies 
T_t\equiv1\pmod{\mathfrak p^r}.
\]
Note that an isometry of $L$ that is congruent to the identity modulo $\mathfrak p^r$ remains so on $L_{n,0}$. Thus, if $t\in\mathfrak p^rL$ then $T_t \in C_n(r)$ (cf. Remark \ref{rmk-para-conj}).

Now let $\bar\ell$ be an isotropic line in $L/\mathfrak p^rL$, and choose a generator $\bar v$. Since $\bar\ell$ is a rank one direct summand, $\bar v$ is primitive (i.e., its image in $L/\mathfrak pL$ is nonzero). Write
\[
\bar v=\sum_{i=1}^{n-1}
(\bar a_i e_i+\bar b_i f_i)+\bar w,
\]
with $\bar a_i,\bar b_i\in\mathfrak o/\mathfrak p^r$, 
where $\bar w$ lies in the remaining one (or two) dimensional summand.

We claim that at least one of the hyperbolic coordinates $\bar a_i,\bar b_i$ is a unit in $\mathfrak o/\mathfrak p^r$. Otherwise, all these coordinates belong to its maximal ideal $\mathfrak p/\mathfrak p^r$, so they vanish upon further reduction modulo $\mathfrak p$. The image of $\bar v$ in $L/\mathfrak pL$ would then lie entirely in the remaining summand. This image is nonzero by primitivity and isotropic because $\bar\ell$  is given to be isotropic. This is a contradiction since we have observed in \S \ref{section-ql} that the remaining summand is anisotropic modulo $\mathfrak p$. 

By permuting hyperbolic pairs and, if necessary, interchanging $e_i$ and $f_i$, we may arrange that the $e_1$-coordinate is a unit. These operations can be taken in $SO(L)$; an interchange of $e_i$ and $f_i$ has determinant $-1$, which we can compensate by a sign change on an orthogonal unit norm basis vector in the remaining summand.

Thus, after multiplying by a unit if necessary, we have 
\[\bar{v} = e_1+\bar t+b f_1\]
with 
$\bar t\in
\bigl(L\cap\langle e_1,f_1\rangle^\perp\bigr)/
\mathfrak p^r\bigl(L\cap\langle e_1,f_1\rangle^\perp\bigr).
$
Since $0=Q(\bar{v})=b+Q(\bar{t})$, we get $b = -Q(\bar{t})$. We choose an integral lift $t$ of $\bar {t}$ so that 
\[
T_t(e_1)=e_1+t-Q(t)f_1
\]
is primitive and isotropic. Its reduction generates the given line $\bar{\ell} \in L/\mathfrak p^rL$ after our coordinate changes. The claim on transitivity amounts to undoing these changes. The same argument over $\mathfrak o$ proves that $C$ and $C\cap H_n(\mathfrak o)$ act transitively on primitive integral isotropic lines in their respective lattices.

Now suppose $\mathfrak o v_1$ and $\mathfrak o v_2$ are two lifts of the same line in $X_r$. By transitivity there exists $h\in C$ carrying $\mathfrak o v_1$ to $\mathfrak o e_1$. The line $h\mathfrak o v_2 \equiv h\mathfrak o v_1 = \mathfrak o e_1 \pmod{\mathfrak p^r}$ thus reduces to $(\mathfrak o/\mathfrak p^r)e_1$. Choose a generator of $h\mathfrak o v_2$ and normalize its $e_1$-coordinate to $1$ as earlier. As in the above argument, we get
\[
h\mathfrak o v_2
=\mathfrak o\bigl(e_1+t-Q(t)f_1\bigr),\]
where $t\in\mathfrak p^r
\bigl(L_{n,0}\cap\xi_{\alpha,u}^{\perp}
\cap\langle e_1,f_1\rangle^\perp\bigr)$.
So we see that $T_t\in C\cap C_n(r)$ carries $\mathfrak o e_1$ to $h\mathfrak o v_2$. Since $C\cap C_n(r)$ is normal in $C$,
\[
h^{-1}T_th\in C\cap C_n(r)
\]
carries $\mathfrak o v_1$ to $\mathfrak o v_2$.

It remains to prove the asserted equality of sets. By Remark~\ref{rmk-para-conj},
\[
\mathcal P_{n,2r}=H_n(\mathfrak o)C_n(r),
\]
so every element of $\mathcal P_{n,2r}$ fixes $v_0$ modulo $\mathfrak p^r$.
Since $p\in\mathcal P_{n,2r}\subset G_n(\mathfrak o)$, the vector $p^{-1}e_1$ is primitive and integral. Moreover, $pv_0\equiv v_0\pmod{\mathfrak p^rL_{n,0}}$, so
\[
B_n(p^{-1}e_1,v_0)
=B_n(e_1,pv_0)
\equiv B_n(e_1,v_0)=0
\mod{\mathfrak p^r}.
\]
Thus the reduction of $p^{-1}\mathfrak o e_1$ is orthogonal to $\bar{v}_0$. Now suppose 
$c=qp\in C\cap \Omega_n\mathcal P_{n,2r}$. Since $q$ stabilizes $\ell_0$,
we get $c^{-1}\ell_0=p^{-1}\ell_0$. Their primitive lattice lines therefore agree:
\[
c^{-1}\mathfrak o e_1=p^{-1}\mathfrak o e_1,
\]
since each is the intersection of this common $F$-line with $L_{n,0}$. Hence the reduction of $c^{-1}\mathfrak o e_1$ is also orthogonal to $\bar v_0$.
Note also that since $c\xi_{\alpha,u}=\xi_{\alpha,u}$, we get 
\[
B_n(c^{-1}e_1,\xi_{\alpha,u})
=B_n(e_1,c\xi_{\alpha,u})
=B_n(e_1,\xi_{\alpha,u})=0.
\]
Here the last equality follows because $n\ge2$ and $\xi_{\alpha,u}\in\langle e_n,v_0,f_n\rangle$. This orthogonality holds over $\mathfrak o$, and hence also modulo $\mathfrak p^r$. We have thus proved that $\rho_r(c) \in Y_r$.

Conversely, suppose that $\rho_r(c)\in Y_r$. Since $\rho_r(1)$ is also in $Y_r$, by transitivity on $Y_r$ we get $d\in C\cap H_n(\mathfrak o)$ such that
\[
\overline{d^{-1}\mathfrak o e_1}
=\overline{c^{-1}\mathfrak o e_1}.
\]
Thus
$d^{-1}\mathfrak o e_1$ and $c^{-1}\mathfrak o e_1$
are two integral lifts of the same line in $Y_r \subset X_r$ and therefore we know that $k\in C\cap C_n(r)$ with
\[
k^{-1}d^{-1}\mathfrak o e_1=c^{-1}\mathfrak o e_1.
\]
Put $p=dk$. Then
$
p\in C\cap\mathcal P_{n,2r}$ and $
p^{-1}\ell_0=c^{-1}\ell_0$.
Hence $q=cp^{-1}$ belongs to $C\cap \Omega_n$, and $c=qp$. We have thus proved
\[
C\cap \Omega_n\mathcal P_{n,2r}
\subseteq\rho_r^{-1}(Y_r)
\subseteq(C\cap \Omega_n)(C\cap\mathcal P_{n,2r}).
\]
The reverse inclusion into the first set is immediate, proving all three sets equal.
\end{proof}

The following cardinality count is a corollary to the first part of Lemma \ref{lemma-dc}.
\begin{corollary}\label{coro-fiber}
We have
\[
|X_r| = q^{(r-1)(2n-2)}|X_1|, ~
|Y_r|=q^{(r-1)(2n-3)}|Y_1|.
\]
\end{corollary}

\begin{proof}
By transitivity established in Lemma \ref{lemma-dc} and equivariance of reduction under $C$, respectively $C \cap H_n(\mathfrak o)$, the fibers of $X_r \twoheadrightarrow X_1$, respectively $Y_r\twoheadrightarrow Y_1$,  have equal cardinality. The normal form established in its proof identifies the fiber over $\mathbb F_q\bar e_1$ with
$\mathfrak pM/\mathfrak p^rM$ via
\[
t\longmapsto \left\langle e_1+t-Q(t)f_1\right\rangle,
\]
where $M = L \cap \langle e_1,f_1\rangle^\perp$. This fiber therefore has cardinality $q^{(r-1) {\rm rank}M}$.
\end{proof}

\begin{remark}\label{kneser}
The material presented in \S \ref{section-ql} and \S \ref{section-lemma-dc} are quite standard. The transitivity and lifting assertions in Lemma \ref{lemma-dc} can also be deduced from integral Witt extension and Kneser’s Hensel lemma; see \cite[Theorem 1.43 and Theorem 5.15]{sch21}, whose proofs follow \cite{kne02}. The assertion on $C \cap C_n(r)$ can be proved as well by following the proof of \cite[Theorem 1.43]{sch21}. Another approach to deduce Corollary \ref{coro-fiber} is via \cite[Lemma 8.26]{sch21}.
\end{remark}

From Corollary \ref{coro-fiber}, we can now compute $|X_r|$ and $|Y_r|$.
\begin{proposition}\label{prop-cardinality}
For $n \geq 2$ and $r \geq 1$,
\[
\frac{|Y_r|}{|X_r|}
=q^{1-r}\frac{q^{n-1}+1}{q^n+1}.
\]
\end{proposition}

\begin{proof}
By Corollary \ref{coro-fiber}, we have
\[
\frac{|Y_r|}{|X_r|}
=q^{1-r}\frac{|Y_1|}{|X_1|}.
\]
The quadratic space defining $X_1$ is nonsplit of dimension $2n$, while that defining $Y_1$ is nondegenerate of dimension $2n-1$. Now applying the quadratic form counts in \cite[Theorems 6.26-6.27]{ln97}, removing the zero vector, and dividing by $q-1$, we get
\[
|X_1|
=\frac{q^{2n-1}-(q-1)q^{n-1}-1}{q-1}
=\frac{(q^{n-1}-1)(q^n+1)}{q-1}, ~
|Y_1|=\frac{q^{2n-2}-1}{q-1}.
\]
\end{proof}

\subsection{Proof of Theorem \ref{thm2}}\label{sub-pf-thm2}

We start with an elementary observation which we state as a lemma whose proof is obvious.
\begin{lemma}\label{lemma-ic}
Let $C$ be a compact group with Haar measure normalized by ${\rm vol}(C)=1$. Suppose $C$ acts transitively on a nonempty finite set $X$, with open point stabilizers. Fix $x_0\in X$. Then, for every function $f:X\to\mathbb C$,
\[
\int_C f(c x_0)dc
=\frac{1}{|X|}\sum_{x\in X}f(x).
\]
\end{lemma}

\begin{proof}[Proof of Theorem \ref{thm2}]
For part (1), Propositions \ref{prop-ggp} and \ref{prop-split-source} give the required representation
$\tau$. Its $K_{1,\alpha}$-distinction follows from
Lemma \ref{lemma-yet} when $\tau$ is tempered. Otherwise, $\tau$ is either an irreducible principal
series ${\rm Ps}(\mu,\mu^{-1})$ or a character $\chi\circ\det$, with $\chi^2=1$. The principal series as
well as the two unramified characters are well-known to be $K_{1,\alpha}$-distinguished. For the
other two ramified characters, it is easy to see that $\tau^{P_{1,m}}=0$ for every $m\geq0$. By \cite[Chapter III, §4, Proposition 4.1(iii)]{kud96}, $\pi_n$ is a subquotient of an
iterated normalized parabolic induction of $\tau$ by unramified characters of ${\rm GL}_1(F)$.
As in the proofs of Propositions \ref{prop-split-source} and \ref{prop-cheng}, successive applications
of Lemma \ref{lemma-cheng} and exactness of compact open invariants then give $\pi_n^{P_{n,m}}
=0$ for every $m\geq0$, which is a contradiction. The converse direction in part (1) follows by Remark \ref{rmk-tempered} and Theorem \ref{thm-cheng}.

For part (2), by Proposition \ref{prop-rsj}, the vectors
\[
v_{j,2r}=q_j(F_{j,2r})
\]
are nonzero paramodular newforms. Put (cf. \S \ref{ortho-para})
\[
h_m=i_m(g_u^1),
~
C_m=h_mC_{m,\alpha}h_m^{-1},
\]
and write $\Lambda_{m,\alpha,u}$ for $\Lambda_{m,\alpha,h_m}$ (cf. (\ref{eq-translated})). Give each $C_m$ Haar measure of total volume one.

We will use the lattice and distinguished vector constructions of \S \ref{section-ql} at each rank $m$, replacing $n$ by $m$, and suppress the rank in the notation $\xi_{\alpha,u}$. Fix $1\leq j <n$. These distinguished vectors $\xi_{\alpha,u}$ correspond under the standard identification
\[
V_j\simeq\langle e_1,f_1\rangle^\perp\subset V_{j+1}.
\]
This will cover both ranks $j+1$ and $j$ throughout the proof.

By the compatibility of the quotient maps (\ref{eq-compatible}), the recursive formula (\ref{eq-recursive}) and the definition of the translated functional (\ref{eq-translated}) give
\[
\Lambda_{j+1,\alpha,u}(v_{j+1,2r})
=
\int_{C_{j+1,\alpha}}
\Lambda_{j,\alpha}\!\left(
q_j\bigl(F_{j+1,2r}(c h_{j+1}^{-1})\bigr)
\right)dc.
\]
We make the change of variables $c^\prime = h_{j+1}ch_{j+1}^{-1}$. Since $h_{j+1}=m(1,h_j)$, from the transformation law (\ref{eq-induction}) and the equivariance of $q_j$, we get 
\[
q_j\bigl(F_{j+1,2r}(h_{j+1}^{-1}c')\bigr)
=
\pi_j(h_j^{-1})
q_j\bigl(F_{j+1,2r}(c')\bigr).
\]
Using (\ref{eq-translated}) again, we get
\[
\Lambda_{j+1,\alpha,u}(v_{j+1,2r})
=
\int_{C_{j+1}}
\Lambda_{j,\alpha,u}\!\left(
q_j\bigl(F_{j+1,2r}(c)\bigr)
\right)dc.
\]

Note that the integrand vanishes outside 
$C_{j+1}\cap\Omega_{j+1}\mathcal P_{j+1,2r}$ by (\ref{eq-paranormal}).
For $c$ in this set, Lemma \ref{lemma-dc} gives
\[
c=pk\]
where
$
p\in C_{j+1}\cap\Omega_{j+1}$ and $k\in C_{j+1}\cap\mathcal P_{j+1,2r}$.
Write the Levi component of $p$ as $m(a,h)$. Since $p$ and $p^{-1}$ preserve the lattice $L_{j+1,0}$ and the line $Fe_1$ they preserve their intersection which is $\mathfrak o e_1$ and thus we have $a\in\mathfrak o^\times$. From $p^{\pm 1}\xi_{\alpha,u}=\xi_{\alpha,u}$, we get $h^{\pm 1} \xi_{\alpha,u}=\xi_{\alpha,u}$ and it follows from the integrality of $p^{\pm 1}$ that $h \in G_j(\mathfrak o)$. Thus,
\[
h \in {\rm Stab}_{G_j(\mathfrak o)}(\xi_{\alpha,u})=C_j.
\]
The normalization (\ref{eq-paranormal}) therefore gives
\[
q_j\bigl(F_{j+1,2r}(c)\bigr)
=
|a|^{2j}
q_j\bigl(\mathcal I_j(h)F_{j,2r}\bigr) = \pi_j(h)v_{j,2r},
\]
as $q_j:\mathcal I_j\to\pi_j$ is a $G_j$-equivariant quotient map.
From the $C_j$-invariance of $\Lambda_{j,\alpha,u}$, we now conclude 
\[
\Lambda_{j,\alpha,u}\!\left(
q_j\bigl(F_{j+1,2r}(c)\bigr)
\right)
=
\Lambda_{j,\alpha,u}(v_{j,2r}).
\]
Since 
\[
\Lambda_{j,\alpha,u}\!\left(q_j(F_{j+1,2r}(c))\right)
=
\begin{cases}
\Lambda_{j,\alpha,u}(v_{j,2r}),&c\in C_{j+1}\cap\Omega_{j+1}\mathcal P_{j+1,2r},\\
0,&c\notin C_{j+1}\cap\Omega_{j+1}\mathcal P_{j+1,2r},
\end{cases}
\]
it follows that
\begin{equation}\label{eq-parallel}
\Lambda_{j+1,\alpha,u}(v_{j+1,2r})
=
\operatorname{vol}\bigl(
C_{j+1}\cap\Omega_{j+1}\mathcal P_{j+1,2r}
\bigr)
\Lambda_{j,\alpha,u}(v_{j,2r}).
\end{equation}

We compute this volume by using the sets $X_r$ and $Y_r$ from \S \ref{section-lemma-dc} now at rank $j+1$. Recall that the map $\rho_r:C \rightarrow X_r$ is defined by $
\rho_r(c)=\overline{c^{-1}\mathfrak o e_1}$.
By Lemma \ref{lemma-dc}, we have transitivity of $C_{j+1}$ on $X_r$ and also
\[
C_{j+1}\cap\Omega_{j+1}\mathcal P_{j+1,2r}
=
\rho_r^{-1}(Y_r).
\]
By transitivity, its fibers are the $|X_r|$ right cosets of the stabilizer of $\overline{\mathfrak oe_1}$. Each stabilizer is open as can be seen as follows. If $K_r$ is the kernel of the action of $C_{j+1}$ on the lattice modulo $\mathfrak p^r$ then it is open as in an integral basis it is defined by the finitely many conditions $c_{ab}-\delta_{ab}\in\mathfrak p^r$, 
and these conditions define an open subset of $C_{j+1}$ since the matrix entry functions are continuous and $\mathfrak p^r$ is open in $F$. Since every element of $K_r$ acts as the identity on the reduced lattice and hence fixes every line in $X_r$, we have $K_r\subseteq {\rm Stab}_{C_{j+1}}(\ell)$ for every $\ell \in X_r$. Now Haar invariance and ${\rm vol}(C_{j+1})=1$ show that each fiber has volume $1/|X_r|$.

Therefore, by Lemma \ref{lemma-ic}, we have
\begin{align*}
{\rm vol}\bigl(
C_{j+1}\cap\Omega_{j+1}\mathcal P_{j+1,2r}
\bigr)
&=\int_{C_{j+1}} 1_{Y_r}\bigl(\rho_r(c)\bigr)dc\\
&=\frac{1}{|X_r|}\sum_{x\in X_r} 1_{Y_r}(x)\\
&=\frac{|Y_r|}{|X_r|}\\
&=q^{1-r}\frac{q^j+1}{q^{j+1}+1},
\end{align*}
where the last step is by Proposition \ref{prop-cardinality}.

Putting all together we have got 
\[
\Lambda_{j+1,\alpha,u}(v_{j+1,2r})
=
q^{1-r}\frac{q^j+1}{q^{j+1}+1}
\Lambda_{j,\alpha,u}(v_{j,2r}).
\]
Iterating for $j=1,\ldots,n-1$, we have
\begin{align*}
\Lambda_{n,\alpha,u}(v_{n,2r})
&=
\left(\prod_{j=1}^{n-1}
q^{1-r}\frac{q^j+1}{q^{j+1}+1}\right)
\Lambda_{1,\alpha,u}(v_{1,2r})\\
&=
q^{(n-1)(1-r)}
\frac{q+1}{q^n+1}\,
\Lambda_{1,\alpha,u}(v_{1,2r}).
\end{align*}

By Remark \ref{rmk-para-conj} and the same lifting argument as in the proof of Corollary \ref{coro-fiber}, for $r\ge1$ we have
\[
[G_j(\mathfrak o):\mathcal P_{j,2r}]
=
q^{2j(r-1)}
\frac{|G_j(\mathbb F_q)|}{|H_j(\mathbb F_q)|}
=
q^{2jr-j}(q^j+1).
\]
Consequently,
\[
D_j(r)=q^{-jr}[G_j(\mathfrak o):\mathcal P_{j,2r}]
=q^{j(r-1)}(q^j+1),
\]
and
\[
q^{(n-1)(1-r)}\frac{q+1}{q^n+1}
=\frac{D_1(r)}{D_n(r)}.
\]
Thus
\[
D_n(r)\Lambda_{n,\alpha,u}(v_{n,2r})
=
D_1(r)\Lambda_{1,\alpha,u}(v_{1,2r}),
\]
as asserted.
\end{proof}

\begin{remark}
The argument leading to (\ref{eq-parallel}) in the proof above parallels the finite field calculation in Lemma \ref{lemma-weil}. In both the arguments, the rank reduction factor depends only on $q$ and the rank, together with the level in the $p$-adic case. All dependence on the representation remains in the lower rank expression.
\end{remark}

\begin{remark}\label{moeglin}
If $\pi_n$ is of Arthur type, genericity suffices in place of temperedness in Theorem \ref{thm2} (2). Indeed, a generic nontempered $\tau$ with $a(\tau)>0$ is an irreducible principal series ${\rm Ps}(\mu,\mu^{-1})$, with $\mu$ nonunitary and ramified. By \cite[Chapter III, §4, Proposition 4.1(iii)]{kud96}, $\mu$ and $\mu^{-1}$ are the only ramified entries in the extended cuspidal support of $\pi_n$. By Mœglin’s description of extended support in Arthur packets \cite[§4.1, Proposition on p. 814]{moe09}, we get
\[
\{\mu,\mu^{-1}\}
=
\{\chi|\cdot|^{-1/2},\chi|\cdot|^{1/2}\}
\]
for a unitary character $\chi$, making ${\rm Ps}(\mu,\mu^{-1})$ reducible.
\end{remark}

\section{Proof of Theorem \ref{thm3}}

For $1\le j\le n$, let $\pi_j$ denote the representation of $SO(2j+1,\mathbb F_q)$ as in the statement of Theorem \ref{thm3}. Let 
\[
E_j=\big\langle G_j(\mathfrak o)v_{j,2}\big\rangle,~
\lambda_{j,u}=\Lambda_{j,\alpha,u}|_{E_j}.
\]
Since $0 \neq v_{j,2}$ is fixed by
$\mathcal P_{j,2}=H_j(\mathfrak o)C_j(1)$, and $C_j(1)$ is normal in $G_j(\mathfrak o)$, the subgroup $C_j(1)$ acts trivially on $E_j$. Thus $E_j$ is a representation of
\[
G_j(\mathfrak o)/C_j(1)\simeq G_j(\mathbb F_q).
\]
Moreover, $v_{j,2}$ is fixed by $H_j(\mathfrak o)$, and hence by $H_j(\mathbb F_q)$. By Proposition \ref{prop-cheng}, we get 
\[
E_j^{H_j(\mathbb F_q)}
=E_j^{\mathcal P_{j,2}}
=\mathbb C v_{j,2}.
\]

Note that the inflated representation $\pi_1$ occurs in $\tau|_{G_1(\mathfrak o)}$. Since $\pi_1^{H_1(\mathbb F_q)}\ne0$, this copy contains a nonzero $\mathcal P_{1,2}$-fixed vector. By Proposition \ref{prop-cheng}, we can identify that vector, up to scalar, with $v_{1,2}$; irreducibility then gives $E_1\simeq\pi_1$. 

With the identification $E_1\simeq\pi_1$, we have 
\[
\dim\pi_1^{H_1(\mathbb F_q)}=1,
\]
and then we also have
\[
\dim\pi_1^{K_{1,\alpha}(\mathbb F_q)}=1,
\]
because for ${\rm PGL}_2(\mathbb F_q)$, an irreducible representation with a one-dimensional split-torus fixed space also has a one-dimensional nonsplit torus fixed space as can be verified easily. It follows that $\lambda_{1,0} \not\equiv 0$. 

By Lemma \ref{kriz2}, we have a $\sigma_0\in\mathcal R$ such that
\[
E_1\simeq \pi_1 \simeq \Theta_{1,\psi}(\sigma_0).
\]

At this point, we observe a consequence of the rank one normalization in Theorem \ref{thm3}. Choose a $G_1(\mathbb F_q)$-invariant Hermitian inner product on $E_1$. Since $\lambda_{1,u}$ are the prescribed translates of one fixed functional by elements of $G_1(\mathfrak o)$, the quantity
\[
d_1=\|v_{1,2}\|^2\|\lambda_{1,u}\|^2
\]
is independent of $u$. It follows from the definition of the correlation coefficient that
\[
|\lambda_{1,u}(v_{1,2})|^2
=d_1\,c(\pi_1;H_1,K_{1,\alpha,u}).
\]
Now there is at least one $u$ for which the left hand side is nonzero. Otherwise, 
as the $g_u^1$, $u \in \mathbb F_q$, form a set of double coset representatives for $H_1\backslash G_1/ K_{1,\alpha}$, we get $\lambda_{1,0}(gv_{1,2})=0$ for all $g \in G_1$ which gives $\lambda_{1,0} \equiv 0$ as $\pi_1$ is irreducible. Thus the rank one normalization in Theorem \ref{thm3} implies that $d_1=1$.

Our next task is to identify $E_j$. 
We now define
\[
\pi_j=\Theta_{j,\psi}(\sigma_0)
\]
and we will prove by induction that $E_j \simeq \pi_j$. 
For the ease of notation, from now on we suppress $\mathbb F_q$ and write $K_j=K_{j,\alpha,u}$. 
Note that $\pi_j$ has fixed vectors for both $H_j$ and $K_j$ by Lemma \ref{kriz2}. 
Let
\[
I_j
={\rm Ind}_{\Omega_j}^{G_j}
( 1\boxtimes\pi_{j-1}).
\]
The two orbit Mackey calculation as in the $p$-adic case in \S \ref{mackey} applies in our situation too.
For either $H_j$ or $K_j$, the closed orbit contributes the corresponding lower rank fixed space. The open orbit contributes $\pi_{j-1}^{G_{j-1}} = 0$, since $\pi_{j-1}$ is nontrivial. Thus we have
\[
I_j^{H_j}
\simeq\pi_{j-1}^{H_{j-1}},
~
I_j^{K_j}
\simeq\pi_{j-1}^{K_{j-1}}.
\]
In particular, both the fixed spaces are one-dimensional.

By Lemma \ref{lemma-gt-finite}, we know $\pi_j \hookrightarrow I_j$. Since
\[
\dim I_j^{H_j}
=\dim\pi_j^{H_j}=1,
\]
the $H_j$-fixed line of $I_j$ lies in this copy of $\pi_j$.

If $I_j^+$ denotes the corresponding $p$-adic induced representation it can be checked,
using the inductive identification $E_{j-1}\simeq\pi_{j-1}$, that we get a $G_j(\mathfrak o)$-equivariant embedding
\[
I_j\hookrightarrow\mathcal (I_j^+)^{C_j(1)}.
\]
Its $H_j$-fixed section is precisely the section used in constructing the newvector, and the local quotient map sends it to $v_{j,2}\ne0$. We conclude that the quotient map is therefore nonzero, hence injective, on the irreducible finite representation $\pi_j$. Note that its image is exactly
\[
\langle G_j(\mathfrak o)v_{j,2}\rangle=E_j.
\]

Let
\[
A_j
=\frac{|\Omega_jH_j|}{|G_j|}
=\frac{[H_j:\Omega_j\cap H_j]}{[G_j:\Omega_j]}
=\frac{q^{j-1}+1}{q^j+1},
\]
and
\[
B_j
=\frac{|\Omega_jK_{j,\alpha,u}|}{|G_j|}
=\frac{[K_j:\Omega_j\cap K_j]}
{[G_j:\Omega_j]}
=\frac{q^{j-1}-1}{q^j-1}.
\]

We equip $I_j$ with the normalized counting inner product and thus
\[
\|f\|_{I_j}^2
=\frac{1}{|G_j(\mathbb F_q)|}
\sum_{g\in G_j(\mathbb F_q)}\|f(g)\|_{E_{j-1}}^2.
\]
Consider the restricted inner product on $\pi_j\subset I_j$ and we transport it to $E_j$ through the isomorphism $p_j:\pi_j\to E_j$; thus $p_j$ is an isometry by construction. 

Let $f_j\in\pi_j\subset I_j$ be supported on $\Omega_j H_j$ such that $f_j(1)=v_{j-1,2}$. For $\omega\in\Omega_j$ and $h\in H_j$,
\[
f_j(\omega h)
=(1\boxtimes\pi_{j-1})(\omega)v_{j-1,2}.
\]
Since the inducing representation is unitary,
\[
\|f_j(\omega h)\|^2=\|v_{j-1,2}\|^2.
\]
Thus,
\[
\|v_{j,2}\|^2 =
\|f_j\|_{I_j}^2
=
\frac{|\Omega_j H_j|}
{|G_j|}
\|v_{j-1,2}\|^2
=A_j\|v_{j-1,2}\|^2.
\]

The recursive local period identity reduces on $I_j$ to
\[
L_{j,u}(F)
=\frac{1}{|K_j|}
\sum_{x\in K_j}\lambda_{j-1,u}(F(x)),
\]
whose restriction to the $\pi_j$-summand is $\lambda_{j,u}$. Let $z_{j-1,u}$ realize the linear form $\lambda_{j-1,u}$, and let $Z_{j,u}$ be the $K_j$-fixed section supported on $\Omega_jK_j$, with value $z_{j-1,u}$ at $1$. Then
\[
\|Z_{j,u}\|^2=B_j\|\lambda_{j-1,u}\|^2,
~
L_{j,u}(F)=B_j^{-1}\langle Z_{j,u}, F \rangle.
\]
Since $Z_{j,u}$ belongs to the $\pi_j$-summand, restricting $L_{j,u}$ to this summand preserves its norm. Therefore
\[
\|\lambda_{j,u}\|^2
=B_j^{-1}\|\lambda_{j-1,u}\|^2.
\]
Using $d_1=1$, we obtain
\[
\begin{aligned}
\|v_{n,2}\|^2\|\lambda_{n,u}\|^2
&=\prod_{j=2}^{n}\frac{A_j}{B_j}\\
&=\prod_{j=2}^{n}
\frac{(q^{j-1}+1)(q^j-1)}
     {(q^j+1)(q^{j-1}-1)}\\
&=\frac{(q+1)(q^n-1)}{(q-1)(q^n+1)}.
\end{aligned}
\]
Since $v_{n,2}$ spans the $H_n$-fixed line in $\pi_n$ and $z_{n,u}$ its $K_{n,\alpha,u}$-fixed line, the definition of the correlation coefficient gives
\begin{align*}
|\Lambda_{n,\alpha,u}(v_{n,2})|^2
&=\|v_{n,2}\|^2\|\lambda_{n,u}\|^2
\,c(\pi_n;H_n,K_{n,\alpha,u})\\
&=\frac{(q+1)(q^n-1)}{(q-1)(q^n+1)}
\,c(\pi_n;H_n,K_{n,\alpha,u}).
\end{align*}
By Proposition \ref{kriz},
\[
\frac{\dim\pi_n}{\dim\pi_1}
=\frac{q^{2n}-1}{q^2-1}.
\]
Since $D_j(1)=q^j+1$, we obtain
\[
\frac{(q+1)(q^n-1)}{(q-1)(q^n+1)}
=
\frac{D_1(1)^2}{D_n(1)^2}
\frac{\dim\pi_n}{\dim\pi_1}.
\]
Thus,
\[
D_n(1)^2\left|\Lambda_{n,\alpha,u}(v_{n,2})\right|^2
=
D_1(1)^2\frac{\dim\pi_n}{\dim\pi_1}
c(\pi_n;H_n,K_{n,\alpha,u}),
\]
as required.


\end{document}